\documentclass[reqno,12pt,a4paper]{amsart}

\usepackage{mathrsfs}
\usepackage{amssymb}
\usepackage{amsfonts}
\usepackage{latexsym}
\usepackage{amsthm}
\usepackage{graphicx}
\usepackage{xcolor}

\theoremstyle{plain}
\newtheorem{theorem}{Theorem}[section]
\newtheorem{lemma}[theorem]{Lemma}
\newtheorem{proposition}[theorem]{Proposition}

\theoremstyle{definition}

\newtheorem{remark}[theorem]{Remark}

\newcommand{\lb}{\label}
\newcommand{\er}{\eqref}

\newcommand{\qqq}{\qquad}
\newcommand{\qq}{\quad}
\newcommand{\no}{\noindent}

\newcommand{\s}{\sigma}
\newcommand{\n}{\nu}
\newcommand{\vk}{\varkappa}
\newcommand{\ve}{\varepsilon}
\newcommand{\G}{\Gamma}
\newcommand{\D}{\Delta}

\newcommand{\cA}{{\mathcal A}}
\newcommand{\cB}{{\mathcal B}}
\newcommand{\cC}{{\mathcal C}}
\newcommand{\cE}{{\mathcal E}}
\newcommand{\cF}{{\mathcal F}}
\newcommand{\cG}{{\mathcal G}}
\newcommand{\cI}{{\mathcal I}}

\newcommand{\cS}{{\mathcal S}}

\newcommand{\cV}{{\mathcal V}}

\newcommand{\A}{{\mathbb A}}
\newcommand{\C}{{\mathbb C}}
\newcommand{\N}{{\mathbb N}}
\newcommand{\R}{{\mathbb R}}
\newcommand{\T}{{\mathbb T}}
\newcommand{\Z}{{\mathbb Z}}

\newcommand{\mH}{{\mathscr H}}

\newcommand{\bG}{{\textbf G}}
\newcommand{\bK}{{\textbf K}}
\newcommand{\be}{{\textbf e}}
\newcommand{\bc}{{\textbf c}}

\newcommand{\bx}{{\textbf x}}
\newcommand{\by}{{\textbf y}}

\newcommand{\ga}{\mathfrak{a}}

\newcommand{\gp}{\mathfrak{p}}
\newcommand{\gs}{\mathfrak{s}}

\newcommand{\mm}{\mathrm m}

\newcommand{\bu}{\bullet}
\newcommand{\ts}{\times}
\newcommand{\sm}{\setminus}
\newcommand{\os}{\oplus}
\newcommand{\ol}{\overline}
\newcommand{\wt}{\widetilde}

\newcommand{\lan}{\langle}
\newcommand{\ran}{\rangle}

\newcommand{\iy}{\infty}
\newcommand{\pa}{\partial}

\newcommand{\Iso}{\mathop{\mathrm{Iso}}\nolimits}

\renewcommand{\[}{\begin{equation}}
\renewcommand{\]}{\end{equation}}

\def\l{\lambda}

\def\t{\tau}
\def\o{\omega}
\def\ss{\subset}

\begin{document}
\makeatletter
\@namedef{subjclassname@2020}{\textup{2020} Mathematics Subject Classification}
\makeatother

\title[On isospectral potentials on periodic discrete graphs]{On isospectral potentials on periodic discrete graphs}

\date{\today}
\author[Natalia Saburova]{Natalia Saburova}
\address{Northern (Arctic) Federal University, Severnaya Dvina Emb. 17, Arkhangelsk, 163002, Russia, \ n.saburova@gmail.com, \ n.saburova@narfu.ru}


\begin{abstract}
We consider discrete Schr\"odinger operators with periodic (possibly complex valued) potentials $Q$ on periodic graphs. Using a system of spectral invariants obtained in \cite{S26}, we investigate the sets of isospectral potentials on periodic graphs. In particular, we show that the Floquet spectrum of the Schr\"odinger operators on the hexagonal lattice and on the Kagome lattice generically (with respect to perturbations of $Q$) determines the potential $Q$ uniquely up to the obvious symmetries (translations and rotations) in the periodic graph.
\end{abstract}

\maketitle

\section{Introduction}
\setcounter{equation}{0}
Discrete Schr\"odinger operators on graphs are used in solid state physics to describe electron motion in a crystal in the tight binding approximation \cite{AM76}. The crystalline structure is modeled by a periodic graph.  Vertices represent locations of atoms in a solid, and edges become chemical bonding of atoms. Discrete operators act on the space of square summable functions defined on the graph vertices. The spectrum of the discrete Schr\"odinger operator with a real periodic potential on a periodic graph consists of a finite number of segments of the real axis (spectral bands) separated by gaps. Spectral bands correspond to the range of energy levels that electrons may have. Studying the band structure of the spectrum allows to explain many physical properties of crystals.

Discrete operators on periodic graphs are intensively studied by mathematicians all over the world during few last decades and there are quite a lot of results for the so-called \emph{forward} spectral problem, i.e., studying spectral properties of the operators, see \cite{SS25} and references therein. At the same time \emph{inverse} spectral problems (i.e., recovering operators from some spectral data) for discrete periodic operators are restricted to the case of the simplest periodic graph, the $d$-dimensional lattice $\mathbb{Z}^d$. However, there are many other models of physical interest (for example, the hexagonal lattice of graphene, the Lieb lattice, the Kagome lattice and so on).

Using Floquet theory the spectral analysis of a $\G$-periodic operator $H$ on a graph $\cG$, where $\G$ is a lattice of rank $d$ in $\R^d$, is reduced to studying spectra of a family of \emph{Floquet} operators $H(k)$ acting on the quotient graph $\cG_*=\cG/\G$ and depending on the \emph{quasimomentum} $k\in\T^d:=(-\pi,\pi]^d$. The family of the spectra $\s(H(k))$, $k\in\T^d$, is called the \emph{Floquet spectrum} of the operator $H$, and the spectrum of $H(0)$ is called the \emph{periodic spectrum} of $H$.

One of interesting and important inverse spectral problems is to understand when two potentials are \emph{isospectral}, i.e., the corresponding Schr\"odinger operators have the same periodic or even Floquet spectrum. This problem has a rich history of study in both discrete and continuous settings, see the surveys \cite{K16,L22}. Isospectral periodic potentials on the $d$-dimensional lattice $\mathbb{Z}^d$, $d\geq2$, were studied by T.Kappeler in a series of papers \cite{Ka88a,Ka88b,Ka89} and by W.Liu \cite{L23} (see also the book \cite{GKT93}). To the best of our knowledge there are no such kind of results for other periodic graphs (except for the triangular lattice also considered in \cite{L23}).

In this paper, we consider the discrete Schr\"odinger operator with periodic (possibly complex valued) potentials on periodic graphs. The goal of this paper is twofold. Firstly, it is an attempt to generalize the results obtained by T.Kappeler \cite{Ka88a,Ka88b,Ka89} for the $d$-dimensional lattice $\Z^d$ to other periodic graphs. Recently, in \cite{S26} we have constructed a system of Floquet spectral invariants for the discrete Schr\"odinger operators on periodic graphs. These invariants are polynomials in the potential and determined by cycles in the quotient graph from some specific cycle sets. The second goal of the present paper is to show how this system of spectral invariants can be used to study the sets of isospectral periodic potentials for the Schr\"odinger operators on periodic graphs.

\medskip

The paper is organized as follows. In Introduction we define the discrete Schr\"odinger operator on periodic graphs and briefly describe its spectrum. In Section \ref{Sec2} we formulate our main results:

$\bu$ we describe cardinality of the sets of isospectral potentials on arbitrary periodic graphs (Proposition \ref{Pcar}) and present asymptotics for large isospectral potentials (Proposition \ref{PAsy});

$\bu$ we show that the Floquet spectrum of the Schr\"odinger operators on the hexagonal lattice and on the Kagome lattice generically (with respect to perturbations of $Q$) determines the potential $Q$ uniquely up to the obvious symmetries (translations and rotations) in the periodic graphs (Theorems \ref{THL} and \ref{TKL}).

In Section \ref{Sec3} we recall the spectral invariants for the discrete Schr\"odinger operator from \cite{S26} and some of their properties (Propositions \ref{PSpIn} and \ref{Pfise}). Section \ref{Sec4} is devoted to the proof of the main results.

\subsection{Periodic graphs.}
Let $\cG=(\cV,\cE)$ be a simple (i.e., without loops and multiple edges) connected infinite graph embedded into the space $\R^d$, where $\cV$ is the set of its vertices and $\cE$ is the set of its unoriented edges. Considering each edge in $\cE$ to have two orientations, we introduce the set $\cA$ of all oriented edges. An edge starting at a vertex $u$ and ending at a vertex $v$ from $\cV$ will be denoted as the ordered pair $(u,v)\in\cA$. Let $\ol\be=(v,u)$ be the inverse edge of $\be=(u,v)\in\cA$. We define the degree ${\vk}_v$ of the vertex $v\in\cV$ as the number of all edges from $\cA$ starting at $v$.

Let $\G$ be a lattice of rank $d$ in $\R^d$ with a basis $\{\ga_1,\ldots,\ga_d\}$, i.e.,
$$
\textstyle\G=\Big\{\ga\in\R^d\mid \ga=\sum\limits_{s=1}^dn_s\ga_s, \; (n_s)_{s=1}^d\in\Z^d\Big\},
$$
and let
\[\lb{fuce}
\textstyle\Omega=\Big\{\bx\in\R^d \mid \bx=\sum\limits_{s=1}^dx_s\ga_s, \; (x_s)_{s=1}^d\in
[0,1)^d\Big\}
\]
be the \emph{fundamental cell} of the lattice $\G$. We define the equivalence relation on $\R^d$:
$$
\bx\equiv \by \; (\hspace{-4mm}\mod \G) \qq\Leftrightarrow\qq \bx-\by\in\G \qqq
\forall\, \bx,\by\in\R^d.
$$

We consider \emph{locally finite $\G$-periodic graphs} $\cG$, i.e., graphs satisfying the following conditions:
\begin{itemize}
  \item[1)] $\cG=\cG+\ga$ for any $\ga\in\G$, i.e., $\cG$ is invariant under translation by any vector $\ga\in\G$;
  \item[2)] the quotient graph  $\cG_*=\cG/\G$ is finite.
\end{itemize}
The basis vectors $\ga_1,\ldots,\ga_d$ of the lattice $\G$ are called the {\it periods}  of $\cG$. The quotient graph $\cG_*=\cG/\G$ is a graph on the $d$-dimensional torus $\R^d/\G$. The graph $\cG_*=(\cV_*,\cE_*)$ has the vertex set $\cV_*=\cV/\G$, the set $\cE_*=\cE/\G$ of unoriented edges and the set $\cA_*=\cA/\G$ of doubled oriented edges.

\begin{remark}
A periodic graph $\cG$ could be defined as an abstract infinite graph  equipped with an action of a finitely generated free abelian group (a lattice) $\G$  (see, e.g., Chapter 4 in \cite{BK13}). The embedding of $\cG$ into $\R^d$ is just a simple geometric realization of periodic graphs.
\end{remark}

\subsection{Schr\"odinger operators on periodic graphs.}
Let $\ell^2(\cV)$ be the Hilbert space of all square summable
functions  $f:\cV\to \C$ equipped with the norm
$$
\textstyle\|f\|^2_{\ell^2(\cV)}=\sum\limits_{v\in\cV}|f(v)|^2<\infty.
$$

We consider the discrete Schr\"odinger operator $H$ acting on $\ell^2(\cV)$ and given by
\[
\lb{Sh} H=A+Q,
\]
where $A$ is the \emph{adjacency} operator having the form
\[\lb{ALO}
(A f)(v)=\sum_{(v,u)\in\cA}f(u), \qqq f\in\ell^2(\cV), \qqq v\in\cV,
\]
and $Q$ is a $\G$-periodic potential (possibly complex valued), i.e., it satisfies
\[\lb{ppot}
(Qf)(v)=Q(v)f(v), \qqq Q(v+\ga)=Q(v), \qqq \forall\,(v,\ga)\in\cV\ts\G.
\]
The sum in \er{ALO} is taken over all edges from $\cA$ starting at the
vertex $v$.

\begin{remark}\lb{Re12}
\emph{i}) It is more usual and physically motivated to consider the Schr\"odinger operator $H=\D+Q$, where $\D$ is the discrete \emph{combinatorial Laplacian} given by $\D=\vk-A$. Here $\vk$ is the degree potential, i.e.,
$$
(\vk f)(v)=\vk_v f(v), \qqq f\in\ell^2(\cV), \qqq v\in\cV,
$$
where $\vk_v$ is the degree of the vertex $v$. But since we do not impose any restrictions on the potential $Q$ (except $\G$-periodicity), we may absorb the degree potential $\vk$ into $Q$.

\emph{ii}) If $\cG$ is a $\G$-periodic graph and $Q$ is a $\G'$-periodic  potential, where $\G'$ is a sublattice of $\G$, then instead of the quotient graph $\cG_*=\cG/\G$ of $\cG$ we will deal with its \emph{expanded} quotient graph $\cG_*=\cG/\G'$.
\end{remark}

\subsection{Spectra of Schr\"odinger operators.}
We briefly describe the spectrum of the Schr\"odinger operator $H$ on periodic graphs (for more details, see, e.g., \cite{HN09} or \cite{KS14}). We introduce the Hilbert space
$$
\mH=L^2\Big(\T^{d},\frac{dk}{(2\pi)^d}\,;\ell^2(\cV_*)\Big)
=\int_{\T^{d}}^{\os}\ell^2(\cV_*)\,\frac{dk}{(2\pi)^d}\,, \qqq
\T^d=\R^d/(2\pi\Z)^d,
$$
equipped with the norm
$$
\|g\|^2_{\mH}=\int_{\T^d}\|g(k)\|_{\ell^2(\cV_*)}^2\frac{dk}{(2\pi)^d}\,,\qqq g\in\mH.
$$

The Schr\"odinger operator $H$ on $\ell^2(\cV)$ has the following decomposition into a constant fiber direct integral, see \cite[Theorem 1.1]{KS14},
$$
UH U^{-1}=\int^\oplus_{\T^d}H(k)\frac{dk}{(2\pi)^d}\,,
$$
where $U:\ell^2(\cV)\to\mH$ is some unitary operator (the Gelfand
transform), the parameter $k$ is called the \emph{quasimomentum}.
For each $k\in\T^d$ the \emph{fiber} (or \emph{Floquet}) operator $H(k)$ on $\ell^2(\cV_*)$ is given by
$$
H(k)=A(k)+Q.
$$
Here $Q$ is the potential on $\ell^2(\cV_*)$, and $A(k)$ is the fiber adjacency operator having the form
\[
\label{fado}
\big(A(k)f\big)(v)=\sum_{\be=(v,u)\in\cA_*}e^{i\lan\t(\be),\,k\ran}f(u),
 \qqq f\in\ell^2(\cV_*),\qqq v\in \cV_*,
\]
where $\t(\be)\in\Z^d$ is the \emph{index} of the edge $\be\in\cA_*$ defined in Section \ref{Sec2.1}, and $\lan\cdot\,,\cdot\ran$ denotes the standard inner product in $\R^d$.

Let $\#M$ denote the number of elements in a set $M$. For each $k\in\T^d$ the spectrum $\s(H(k))$ of $H(k)$ consists of $\n:=\#\cV_*$ eigenvalues (including multiplicities). The family of spectra of all Floquet operators $H(k)$
$$
\big\{\s\big(H(k)\big), k\in\T^d\big\}
$$
is called the \emph{Floquet} spectrum of the Schr\"odinger operator $H$. The Floquet operator $H(0)$ is just the Schr\"odinger operator defined by \er{Sh} -- \er{ppot} on the finite quotient graph $\cG_*$. Its spectrum $\s\big(H(0)\big)$ is called the \emph{periodic} spectrum of $H$.

If the potential $Q$ is real, then each fiber operator $H(k)$, $k\in\T^{d}$, is self-adjoint and has $\n$ real eigenvalues $\l_j(k)$, $j=1,\ldots,\n$, labeled in non-decreasing order counting multiplicities:
$$
\l_{1}(k)\leq\l_{2}(k)\leq\ldots\leq\l_{\nu}(k), \qqq
\forall\,k\in\T^{d}.
$$
Each \emph{band function} $\l_j(\cdot)$ is a continuous
and piecewise real analytic function on the torus $\T^{d}$ and creates the
\emph{spectral band} $\s_j(H)$ given by
\[\lb{ban.1H}
\begin{array}{l}
\s_j(H)=[\l_j^-,\l_j^+]=\l_j(\T^{d}),\qqq j\in\N_\n=\{1,\ldots,\n\},\\[8pt]
\displaystyle\textrm{where}\qqq \l_j^-=\min_{k\in\T^d}\l_j(k),\qqq \l_j^+=\max_{k\in\T^d}\l_j(k).
\end{array}
\]
The spectrum of the Schr\"odinger operator $H$ with a real periodic potential $Q$ on a periodic graph $\cG$ has the form
$$
\s(H)=\bigcup_{k\in\T^d}\s\big(H(k)\big)=
\bigcup_{j=1}^{\nu}\s_j(H),
$$
i.e., it consists of $\n$ bands $\s_j(H)$ defined by \er{ban.1H}.

\begin{remark}
Some of spectral bands $\s_j(H)=[\l_j^-,\l_j^+]$ may be degenerate, i.e., $\l_j^-=\l_j^+$.
\end{remark}

\section{Main results}\lb{Sec2}
\setcounter{equation}{0}

\subsection{Isospectral potentials} Let $H=A+Q$ be the Schr\"odinger operator with a $\G$-periodic possibly complex-valued potential $Q$ on a $\G$-periodic graph $\cG$.

Two $\G$-periodic potentials $Q_1$ and $Q_2$ on the same periodic graph $\cG$ are called

\no $\bu$ \emph{Floquet isospectral}, if $A+Q_1$ and $A+Q_2$ have the same Floquet spectrum, i.e.,
$$
\s\big(A(k)+Q_1\big)=\s\big(A(k)+Q_2\big) \qqq \forall\,k\in\T^d;
$$

\no $\bu$ simply \emph{isospectral}, if $A+Q_1$ and $A+Q_2$ have the same periodic spectrum, i.e.,
$$
\s\big(A(0)+Q_1\big)=\s\big(A(0)+Q_2\big),
$$
where $A(k)$ is the Floquet operator defined by \er{fado}.

For a periodic potential $Q$, we denote by $\Iso_\R(Q)$ and  $\Iso_\C(Q)$ (respectively, by $\Iso_\R^{\mathrm{F}}(Q)$ and $\Iso_\C^{\mathrm{F}}(Q)$) the sets of all \emph{real} and \emph{complex} potentials isospectral (respectively, Floquet isospectral) to $Q$. The following proposition describes cardinality of the sets of isospectral potentials on a periodic graph $\cG$.

\begin{proposition}\lb{Pcar} Let $H=A+Q$ be the Schr\"odinger operator with a $\G$-periodic potential $Q$ on a periodic graph $\cG$ with the quotient graph $\cG_*=\cG/\G=(\cV_*,\cE_*)$, $\cV_*=\{v_1,\ldots,v_\nu\}$. Then \\[4pt]
i) $\#\Iso_\C(Q)\leq\nu!$, where $\nu=\#\mathcal{V}_*$;\\[4pt]
ii) $\#\Iso_\C(Q)=\nu!$ generically; \\[4pt]
iii) $\#\Iso_\R(Q)=\nu!$ for real $Q$ such that
\[\lb{sepQ}
d(Q):=\min\limits_{i\neq j}|q_i-q_j|>4\rho, \qqq q_j=Q(v_j),\qqq j\in\N_\n,
\]
where $\rho$ is the spectral radius of the adjacency operator $A(0)$ on the quotient graph $\cG_*$. Moreover, each of $\n!$ regions of $\R^\n$ given by
\[\lb{repe}
\R^\n_\pi:=\big\{(x_1,\ldots,x_\n)\in\R^\n \mid x_{\pi(1)}<\ldots<x_{\pi(\nu)}\big\},
\]
where $\pi$ is a permutation of the set $\N_\n$, contains exactly one potential isospectral to $Q$.
\end{proposition}

\begin{remark}
\emph{i}) "Generically" means outside a proper algebraic subset of the potentials space $\C^\nu$. Recall that a $\G$-periodic potential $Q$ is uniquely determined by its values at the vertices of the quotient graph $\cG_*=\cG/\G=(\cV_*,\cE_*)$, $\n=\#\cV_*$.

\emph{ii}) The spectral radius of a square matrix is the largest absolute value of its eigenvalues. Since $\s(A(0))\subseteq[-\vk_o,\vk_o]$, where $\vk_o$ is the maximum vertex degree, we can replace the condition $d(Q)>4\rho$ by the  weaker (but simpler checked) condition $d(Q)>4\vk_o$.

\emph{iii}) The condition $d(Q)>4\rho$ (or $d(Q)>4\vk_o$), in particular, implies that the spectral bands of the Schr\"odinger operator $H=A+Q$ on the periodic graph $\cG$ do not overlap.

\emph{iv}) For the lattice $\Z^d$, almost the same result was proved in \cite{Ka88a}. The only difference is in \emph{iii}). It was shown that $\#\Iso_\R(Q)=\nu!$ for \emph{sufficiently large} real potentials $Q$ with $q_i\neq q_j$, $i\neq j$, \cite[Proposition 2.2]{Ka88a}, without specification how large the potentials should be.

\emph{v}) The regions \er{repe} are the $\nu!$ open connected regions into which $\R^\n$ is divided by the hyperplanes $x_i=x_j$, $i\neq j$.
\end{remark}

The following statement gives asymptotics for large isospectral potentials.

\begin{proposition}\lb{PAsy}
Let $H=A+tQ$ be the Schr\"odinger operator with a real $\G$-periodic potential $Q$ on a periodic graph $\cG$, where $t>0$ is the coupling constant. Assume that the potential $Q$ has pairwise distinct values at the vertices of the quotient graph $\cG_*=\cG/\G=(\cV_*,\cE_*)$, $\cV_*=\{v_1,\ldots,v_\nu\}$, and $\cG_*$ has no loops. Then for $t>t_o:=\frac{4\rho}{d(Q)}$, where
$$
d(Q):=\min\limits_{i\neq j}|q_i-q_j|, \qqq q_j=Q(v_j),\qqq j\in\N_\n,
$$
and $\rho$ is the spectral radius of the adjacency operator $A(0)$ on $\cG_*$, the set of potentials isospectral to $tQ$ consists of $\n!$ real potentials $t P$ and
\[\lb{asy}
\textstyle  p_i=P(v_i)=q_{\pi(i)}+p_i^{(2)}(Q)\,\frac{1}{t^2}+ O\big(\frac1{t^3}\big), \qqq i\in\N_\n,\qqq \textrm{as}\qqq t\to\iy.
\]
Here $\pi$ is a permutation of the set $\N_\n$ (one isospectral potential $t P$ for each of $\n!$ permutations~$\pi$);
\[\lb{coef2}
p_i^{(2)}(Q)=\frac1{\prod\limits_{\substack{k=1\\ k\neq i}}^\n(q_{\pi(i)}-q_{\pi(k)})}\sum_{j=1}^{\n-2}(-1)^{\n-j}
\gs_{\n-j-2}(q_{\pi(1)},\ldots,\hat q_{\pi(i)},\ldots,q_{\pi(\n)})
b_j(Q),
\]
$$
b_j(Q)=\sum_{\be=\{v_l,v_k\}\in\cE_*}m_\be^2
\big(h_j(q_l,q_k)-h_j(q_{\pi(l)},q_{\pi(k)})\big),\qqq
h_j(q_l,q_k)=\sum_{\substack{j_1+j_2=j\\j_1,\,j_2\geq0}}q_l^{j_1}q_k^{j_2},
$$
$$
\gs_k(x_1,\ldots,x_n)=\sum_{1\leq a_1<\ldots<a_k\leq n}x_{a_1}\ldots x_{a_k}
$$
is the elementary symmetric polynomial of degree $k$ in variables $x_1,\ldots,x_n$, and $\hat x_i$ means that the argument $x_i$ is excluded; $m_\be$ is the multiplicity of the edge $\be\in\cE_*$.
\end{proposition}

\begin{remark}
The asymptotics (\ref{asy}) shows that, for a large real potential $Q=(q_1,\ldots,q_\nu)$, each of its isospectral potential $P$ is obtained from $Q$ by a permutation of its components plus a small term depending on the structure of the quotient graph $\cG_*$. In particular, if the bijection $\pi_*:\cV_*\to\cV_*$ defined by $\pi_*(v_i)=v_{\pi(i)}$, $i\in\N_\nu$, is an automorphism of the quotient graph $\cG_*$, then, in \er{coef2}, $b_j(Q)=0$, $j\in\N_{\nu-2}$. Indeed, in this case, the potential
$P=(q_{\pi(1)},\ldots,q_{\pi(\n)})$ is isospectral to $Q$.
\end{remark}

\subsection{Floquet isospectrality on the hexagonal lattice}
In this section we consider the discrete Schr\"odinger operator $H=A+Q$ with periodic potentials $Q$ on the hexagonal lattice $\bG$ and describe the sets of Floquet isospectral potentials on $\bG$.

\begin{figure}[t!]\centering
\unitlength 1.2mm
\begin{picture}(55,45)(0,0)
\put(2,36){$\bG$}

\put(9,3){\circle{1.3}}
\put(21,3){\circle{1.3}}
\put(33,3){\circle{1.3}}
\put(45,3){\circle{1.3}}

\put(3,6){\circle*{1.3}}
\put(15,6){\circle*{1.3}}
\put(27,6){\circle*{1.3}}
\put(39,6){\circle*{1.3}}
\put(51,6){\circle*{1.3}}

\put(3,12){\circle{1.3}}
\put(15,12){\circle{1.3}}
\put(27,12){\circle{1.3}}
\put(39,12){\circle{1.3}}
\put(51,12){\circle{1.3}}

\put(9,15){\circle*{1.3}}
\put(21,15){\circle*{1.3}}
\put(33,15){\circle*{1.3}}
\put(45,15){\circle*{1.3}}
\put(9,21){\circle{1.3}}
\put(21,21){\circle{1.3}}
\put(33,21){\circle{1.3}}
\put(45,21){\circle{1.3}}

\put(3,24){\circle*{1.3}}
\put(15,24){\circle*{1.3}}
\put(27,24){\circle*{1.3}}
\put(39,24){\circle*{1.3}}
\put(51,24){\circle*{1.3}}

\put(3,30){\circle{1.3}}
\put(15,30){\circle{1.3}}
\put(27,30){\circle{1.3}}
\put(39,30){\circle{1.3}}
\put(51,30){\circle{1.3}}

\put(9,33){\circle*{1.3}}
\put(21,33){\circle*{1.3}}
\put(33,33){\circle*{1.3}}
\put(45,33){\circle*{1.3}}
\put(9,39){\circle{1.3}}
\put(21,39){\circle{1.3}}
\put(33,39){\circle{1.3}}
\put(45,39){\circle{1.3}}

\put(3,42){\circle*{1.3}}
\put(15,42){\circle*{1.3}}
\put(27,42){\circle*{1.3}}
\put(39,42){\circle*{1.3}}
\put(51,42){\circle*{1.3}}

\put(9,3){\line(0,-1){3.00}}
\put(21,3){\line(0,-1){3.00}}
\put(33,3){\line(0,-1){3.00}}
\put(45,3){\line(0,-1){3.00}}

\bezier{60}(3,6)(1.5,5.25)(0,4.5)
\bezier{60}(3,12)(1.5,12.75)(0,13.5)
\bezier{60}(3,24)(1.5,23.25)(0,22.5)
\bezier{60}(3,30)(1.5,30.75)(0,31.5)
\bezier{60}(3,42)(1.5,41.25)(0,40.5)

\bezier{60}(51,6)(52.5,5.25)(54,4.5)
\bezier{60}(51,12)(52.5,12.75)(54,13.5)
\bezier{60}(51,24)(52.5,23.25)(54,22.5)
\bezier{60}(51,30)(52.5,30.75)(54,31.5)
\bezier{60}(51,42)(52.5,41.25)(54,40.5)

\put(9,3){\line(-2,1){6.00}}
\put(9,3){\line(2,1){6.00}}
\put(21,3){\line(-2,1){6.00}}
\put(21,3){\line(2,1){6.00}}
\put(33,3){\line(-2,1){6.00}}
\put(33,3){\line(2,1){6.00}}
\put(45,3){\line(-2,1){6.00}}
\put(45,3){\line(2,1){6.00}}
\put(3,6){\line(0,1){6.00}}
\put(15,6){\line(0,1){6.00}}
\put(27,6){\line(0,1){6.00}}
\put(39,6){\line(0,1){6.00}}
\put(51,6){\line(0,1){6.00}}
\put(9,15){\line(-2,-1){6.00}}
\put(9,15){\line(2,-1){6.00}}
\put(21,15){\line(-2,-1){6.00}}
\put(21,15){\line(2,-1){6.00}}
\put(33,15){\line(-2,-1){6.00}}
\put(33,15){\line(2,-1){6.00}}
\put(45,15){\line(-2,-1){6.00}}
\put(45,15){\line(2,-1){6.00}}
\put(9,15){\line(0,1){6.00}}
\put(21,15){\line(0,1){6.00}}
\put(33,15){\line(0,1){6.00}}
\put(45,15){\line(0,1){6.00}}
\put(9,21){\line(-2,1){6.00}}
\put(9,21){\line(2,1){6.00}}
\put(21,21){\line(-2,1){6.00}}
\put(21,21){\line(2,1){6.00}}
\put(33,21){\line(-2,1){6.00}}
\put(33,21){\line(2,1){6.00}}
\put(45,21){\line(-2,1){6.00}}
\put(45,21){\line(2,1){6.00}}
\put(3,24){\line(0,1){6.00}}
\put(15,24){\line(0,1){6.00}}
\put(27,24){\line(0,1){6.00}}
\put(39,24){\line(0,1){6.00}}
\put(51,24){\line(0,1){6.00}}
\put(9,33){\line(-2,-1){6.00}}
\put(9,33){\line(2,-1){6.00}}
\put(21,33){\line(-2,-1){6.00}}
\put(21,33){\line(2,-1){6.00}}
\put(33,33){\line(-2,-1){6.00}}
\put(33,33){\line(2,-1){6.00}}
\put(45,33){\line(-2,-1){6.00}}
\put(45,33){\line(2,-1){6.00}}
\put(9,33){\line(0,1){6.00}}
\put(21,33){\line(0,1){6.00}}
\put(33,33){\line(0,1){6.00}}
\put(45,33){\line(0,1){6.00}}
\put(9,39){\line(-2,1){6.00}}
\put(9,39){\line(2,1){6.00}}
\put(21,39){\line(-2,1){6.00}}
\put(21,39){\line(2,1){6.00}}
\put(33,39){\line(-2,1){6.00}}
\put(33,39){\line(2,1){6.00}}
\put(45,39){\line(-2,1){6.00}}
\put(45,39){\line(2,1){6.00}}
\put(3,42){\line(0,1){3.00}}
\put(15,42){\line(0,1){3.00}}
\put(27,42){\line(0,1){3.00}}
\put(39,42){\line(0,1){3.00}}
\put(51,42){\line(0,1){3.00}}

\put(3.6,42.7){$\scriptstyle v_{2,-2,4}$}
\put(15.6,42.7){$\scriptstyle v_{2,-1,4}$}
\put(27.5,42.7){$\scriptstyle v_{2,0,4}$}
\put(40,42.7){$\scriptstyle v_{2,1,4}$}
\put(52,42.7){$\scriptstyle v_{2,2,4}$}

\put(9.6,37.5){$\scriptstyle v_{1,-1,4}$}
\put(21.8,37.5){$\scriptstyle v_{1,0,4}$}
\put(33.8,37.5){$\scriptstyle v_{1,1,4}$}
\put(45.8,37.5){$\scriptstyle v_{1,2,4}$}

\put(9.6,33.5){$\scriptstyle v_{2,-1,3}$}
\put(21.8,33.5){$\scriptstyle v_{2,0,3}$}
\put(33.8,33.5){$\scriptstyle v_{2,1,3}$}
\put(45.8,33.5){$\scriptstyle v_{2,2,3}$}

\put(3.6,28.3){$\scriptstyle v_{1,-1,3}$}
\put(15.6,28.3){$\scriptstyle v_{1,0,3}$}
\put(28,28.3){$\scriptstyle v_{1,1,3}$}
\put(40,28.3){$\scriptstyle v_{1,2,3}$}
\put(52,28.3){$\scriptstyle v_{1,3,3}$}

\put(3.6,24.7){$\scriptstyle v_{2,-1,2}$}
\put(15.6,24.7){$\scriptstyle v_{2,0,2}$}
\put(28,24.7){$\scriptstyle v_{2,1,2}$}
\put(40,24.7){$\scriptstyle v_{2,2,2}$}
\put(52,24.7){$\scriptstyle v_{2,3,2}$}

\put(2,19.5){$\scriptstyle v_{1,0,2}$}
\put(21.8,19.5){$\scriptstyle v_{1,1,2}$}
\put(33.8,19.5){$\scriptstyle v_{1,2,2}$}
\put(45.8,19.5){$\scriptstyle v_{1,3,2}$}

\put(2,15.5){$\scriptstyle v_{2,0,1}$}
\put(14.6,15.5){$\scriptstyle v_{2,1,1}$}
\put(33.8,15.5){$\scriptstyle v_{2,2,1}$}
\put(45.8,15.5){$\scriptstyle v_{2,3,1}$}

\put(-4,10.3){$\scriptstyle v_{1,0,1}$}
\put(15.5,10.5){$\scriptstyle v_{1,1,1}$}
\put(28,10.3){$\scriptstyle v_{1,2,1}$}
\put(40,10.3){$\scriptstyle v_{1,3,1}$}
\put(52,10.3){$\scriptstyle v_{1,4,1}$}

\put(-4,6.7){$\scriptstyle v_{2,0,0}$}
\put(8,6.7){$\scriptstyle v_{2,1,0}$}
\put(28,6.7){$\scriptstyle v_{2,2,0}$}
\put(40,6.7){$\scriptstyle v_{2,3,0}$}
\put(52,6.7){$\scriptstyle v_{2,4,0}$}

\put(9.8,1.6){$\scriptstyle v_{1,1,0}$}
\put(21.8,1.6){$\scriptstyle v_{1,2,0}$}
\put(33.8,1.6){$\scriptstyle v_{1,3,0}$}
\put(45.8,1.6){$\scriptstyle v_{1,4,0}$}

\put(17,6.5){\small$\ga_1$}
\put(10.4,16.5){\small$\ga_2$}
\put(6.5,8.7){$\scriptstyle O$}
\put(9,9){\vector(1,0){12.0}}
\put(9,9){\vector(2,3){6}}
\bezier{25}(21,9)(24,13.5)(27,18)
\bezier{30}(15,18)(21,18)(27,18)
\end{picture}
\caption{\scriptsize The hexagonal lattice $\bG$ is invariant under translations by the vectors $\ga_1,\ga_2$. There are two types of vertices in $\bG$: the white and black points (with the first index 1 and 2, respectively); the remaining two vertex indices correspond to shifts along the vectors $\ga_1$ and $\ga_2$. Edges connect only white and black vertices.} \label{fig1}
\end{figure}
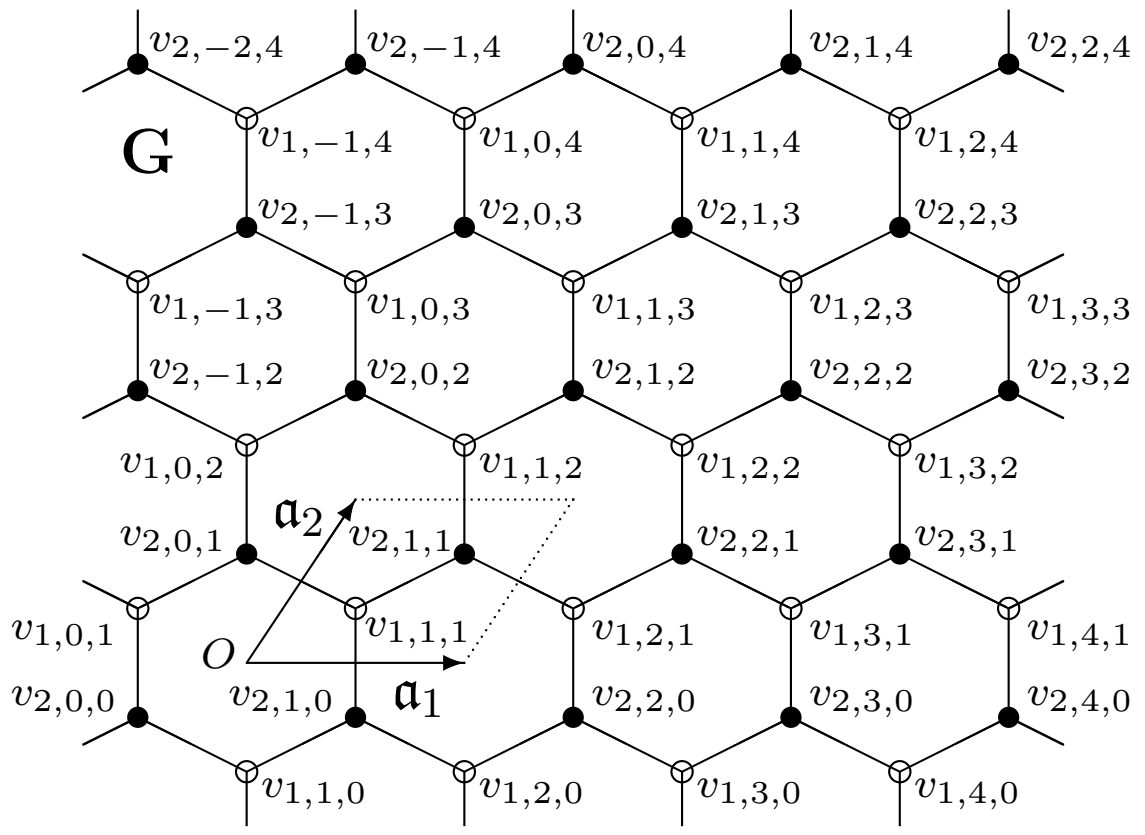

The hexagonal lattice $\bG=(\cV,\cE)$ is a planar periodic graph invariant under translations by the vectors $\ga_1,\ga_2$, see Fig.~\ref{fig1}. The vertex set $\cV$ and the edge set $\cE$ of $\bG$ are given by
$$
\begin{array}{l}
\cV=\{v_{s,n}, \; s=1,2, \;n\in\Z^2\},\\[4pt]
\cE=\big\{(v_{2,n},v_{1,n}),(v_{2,n},v_{1,n_1+1,n_2}),
(v_{2,n},v_{1,n_1,n_2+1}),\qq n=(n_1,n_2)\in\Z^2\big\},
\end{array}
$$
where, with respect to the basis $\ga_1,\ga_2$,
$$
\textstyle v_{1,n}=n-\big(\frac23\,,\frac23\big),\qqq v_{2,n}=n-\big(\frac13\,,\frac13\big),\qqq n\in\Z^2.
$$

Let $H=A+Q$ be the discrete Schr\"odinger operator with a $\G$-periodic potential $Q:\cV\to\C$ on the hexagonal lattice $\bG=(\cV,\cE)$, where $\G$ is the lattice with the basis $p\ga_1,p\ga_2$ for some $p\in\N$. Let
\[\lb{Ga_o}
\G_o=\{\ga=l_1\ga_1+l_2\ga_2 \mid l_j=0,1,\ldots,p-1;\; j=1,2\}.
\]
For a $\G$-periodic potential $Q$ on the hexagonal lattice $\bG$, we define\\
$\bu$ the \emph{translated} potentials $Q_\ga$, $\ga\in\G_o$, by
\[\lb{TrQ}
Q_\ga(v)=Q(v+\ga),
\]
$\bu$ the (half-turn) \emph{rotated} potential $Q^r$ (around the point $O$, see Fig.~\ref{fig1}) by
\[\lb{ReQ}
Q^r(v_{1,n_1,n_2})=Q(v_{2,1-n_1,1-n_2}),\qqq
Q^r(v_{2,n_1,n_2})=Q(v_{1,1-n_1,1-n_2}),
\]
for all $v\in\cV=\{v_{s,n_1,n_2}, \; s=1,2, \;n_1,n_2\in\Z\}$.

\medskip

The following theorem shows that the Floquet spectrum of the Schr\"odinger operator on the hexagonal lattice $\bG$ \emph{generically} (with respect to perturbations of $Q$) determines the periodic potential $Q$ uniquely \emph{up to the obvious symmetries} (translations and half-turn rotation) in~$\bG$.

\begin{theorem}\lb{THL} Let $H=A+Q$ be the Schr\"odinger operator with a $\G$-periodic complex potential $Q$ on the hexagonal lattice $\bG$. Then, generically, the set of all $\G$-periodic potentials Floquet isospectral to $Q$ has the form
\[\lb{FIsG}
\mathrm{Iso}_\C^{\mathrm{F}}(Q)=\big\{Q_\ga,Q_\ga^r \mid \ga\in\G_o\},
\]
where $Q_\ga$ and $Q^r$ are the translated and rotated potentials defined by \er{TrQ} and \er{ReQ}, and $\G_o$ is given by \er{Ga_o}.
\end{theorem}

\begin{remark}
\emph{i}) Even a small portion of the Floquet spectrum does already determine the periodic potential $Q$ generically, again up to the symmetry group of the problem.

\emph{ii}) The Floquet spectrum does only generically determine the potential $Q$ uniquely up to the symmetries. For example, a complex periodic potential  Floquet isospectral to the zero potential on the hexagonal lattice was constructed in \cite{SM25}.
\end{remark}

\subsection{Floquet isospectrality on the Kagome lattice}
\begin{figure}[t!]\centering
\unitlength 1.0mm 
\linethickness{0.4pt}
\begin{picture}(85,58)(0,0)
\put(8.0,51){$\bK$}
\put(-3,15){\line(1,0){66.0}}
\put(7,35){\line(1,0){66.0}}
\put(17,55){\line(1,0){66.0}}
\put(-1.5,12){\line(1,2){23.0}}
\put(18.5,12){\line(1,2){23.0}}
\put(38.5,12){\line(1,2){23.0}}
\put(58.5,12){\line(1,2){23.0}}

\put(11.5,12){\line(-1,2){8.0}}
\put(31.5,12){\line(-1,2){18.0}}
\put(51.5,12){\line(-1,2){23.0}}
\put(66.5,22){\line(-1,2){18.0}}
\put(76.5,42){\line(-1,2){8.0}}

\put(20,15){\vector(1,0){20.0}}
\put(20,15){\vector(1,2){10.0}}
\put(0,15){\circle{1.3}}
\put(10,15){\circle*{1.3}}
\put(20,15){\circle{1.3}}
\put(30,15){\circle*{1.3}}
\put(40,15){\circle{1.3}}
\put(50,15){\circle*{1.3}}
\put(60,15){\circle{1.3}}

\color{blue}
\put(5,25){\circle*{1.3}}
\put(25,25){\circle*{1.3}}
\put(45,25){\circle*{1.3}}
\put(65,25){\circle*{1.3}}

\put(15,45){\circle*{1.3}}
\put(35,45){\circle*{1.3}}
\put(55,45){\circle*{1.3}}
\put(75,45){\circle*{1.3}}
\color{black}
\put(10,35){\circle{1.3}}
\put(20,35){\circle*{1.3}}
\put(30,35){\circle{1.3}}
\put(40,35){\circle*{1.3}}
\put(50,35){\circle{1.3}}
\put(60,35){\circle*{1.3}}
\put(70,35){\circle{1.3}}

\put(20,55){\circle{1.3}}
\put(30,55){\circle*{1.3}}
\put(40,55){\circle{1.3}}
\put(50,55){\circle*{1.3}}
\put(60,55){\circle{1.3}}
\put(70,55){\circle*{1.3}}
\put(80,55){\circle{1.3}}

\put(34.0,12.5){\small$\ga_1$}
\put(24.0,31.4){\small$\ga_2$}
\put(15.9,12.5){$\scriptstyle O$}
\put(20.0,12.5){$\scriptstyle v_{1,1,1}$}
\put(40.0,12.5){$\scriptstyle v_{1,2,1}$}
\put(0,12.5){$\scriptstyle v_{1,0,1}$}
\put(60.0,12.5){$\scriptstyle v_{1,3,1}$}

\put(29.5,16.5){$\scriptstyle v_{2,1,1}$}
\put(49.5,16.5){$\scriptstyle v_{2,2,1}$}
\put(9.5,16.5){$\scriptstyle v_{2,0,1}$}

\put(6.0,25.5){$\scriptstyle v_{3,0,1}$}
\put(26.0,25.5){$\scriptstyle v_{3,1,1}$}
\put(46.0,24.5){$\scriptstyle v_{3,2,1}$}
\put(66.0,25.5){$\scriptstyle v_{3,3,1}$}

\put(39.5,36.0){$\scriptstyle v_{2,1,2}$}
\put(59.5,36.0){$\scriptstyle v_{2,2,2}$}
\put(19.5,36.0){$\scriptstyle v_{2,0,2}$}

\put(9.5,32.5){$\scriptstyle v_{1,0,2}$}
\put(30,32.5){$\scriptstyle v_{1,1,2}$}
\put(50,32.5){$\scriptstyle v_{1,2,2}$}
\put(70,32.5){$\scriptstyle v_{1,3,2}$}

\put(15.9,45.0){$\scriptstyle v_{3,0,2}$}
\put(35.9,45.0){$\scriptstyle v_{3,1,2}$}
\put(55.9,45.0){$\scriptstyle v_{3,2,2}$}
\put(75.9,45.0){$\scriptstyle v_{3,3,2}$}

\put(49.5,56.5){$\scriptstyle v_{2,1,3}$}
\put(69.5,56.5){$\scriptstyle v_{2,2,3}$}
\put(29.5,56.5){$\scriptstyle v_{2,0,3}$}

\put(20,52.5){$\scriptstyle v_{1,0,3}$}
\put(40,52.5){$\scriptstyle v_{1,1,3}$}
\put(60,52.5){$\scriptstyle v_{1,2,3}$}
\put(79.5,52.5){$\scriptstyle v_{1,3,3}$}

\put(20,15){\line(1,0){20.0}}
\put(20,15){\line(1,2){10.0}}
\put(35,45){\line(1,-2){10.0}}
\end{picture}
\vspace{-9mm}
\caption{\scriptsize The Kagome lattice $\bK$; $\ga_1,\ga_2$ are the periods of $\bK$. There are three types of vertices in $\bK$: the white, black and blue points (with the first index 1, 2 and 3, respectively); the remaining two vertex indices correspond to shifts along the vectors $\ga_1$ and $\ga_2$.}
\label{fig2}
\end{figure}
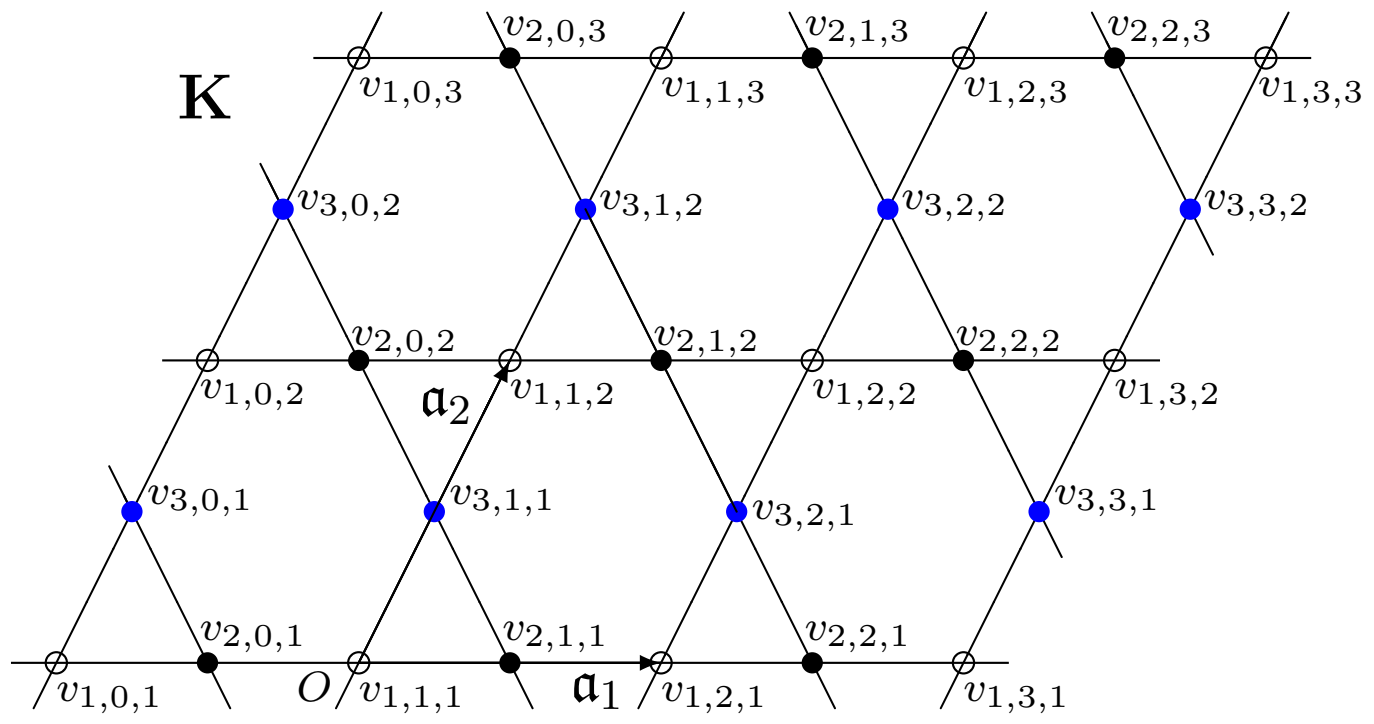

The \emph{Kagome} lattice $\bK=(\cV,\cE)$, named after a traditional Japanese bamboo basket-weaving technique, is a planar periodic graph with periods $\ga_1,\ga_2$, see Fig.~\ref{fig2}. It is known that some minerals contain layers with the Kagome lattice arrangement of atoms in their crystal structure. The vertex set $\cV$ and the edge set $\cE$ of $\bK$ are given by
\begin{multline*}
\begin{array}{l}
\cV=\{v_{s,n}, \; s=1,2,3, \;n\in\Z^2\},\\[6pt]
\cE=\big\{(v_{1,n},v_{2,n}),(v_{1,n},v_{3,n}),
(v_{2,n},v_{3,n}),(v_{2,n},v_{1,n_1+1,n_2}),
\end{array}\\
\hspace{4mm}(v_{3,n},v_{1,n_1,n_2+1}),(v_{3,n},v_{2,n_1-1,n_2+1}),\; n=(n_1,n_2)\in\Z^2\big\},
\end{multline*}
where, with respect to the basis $\ga_1,\ga_2$,
$$
\textstyle v_{1,n}=n-(1,1),\quad v_{2,n}=n-\big(\frac12\,,1\big),\qq v_{3,n}=n-\big(1,\frac12\,\big).
$$

We consider the discrete Schr\"odinger operator $H=A+Q$ with a $\G$-periodic potential $Q:\cV\to\C$ on the Kagome lattice $\bK=(\cV,\cE)$, where $\G$ is the lattice with the basis $p\ga_1,p\ga_2$ for some $p\in\N$. Let
\[\lb{Ga_o'}
\G_o=\{\ga=l_1\ga_1+l_2\ga_2 \mid l_j=0,1,\ldots,p-1;\; j=1,2\}.
\]
For a $\G$-periodic potential $Q$ on the Kagome lattice $\bK$, we define\\
$\bu$ the \emph{translated} potentials $Q_\ga$, $\ga\in\G_o$, by
\[\lb{TrQ'}
Q_\ga(v)=Q(v+\ga),
\]
$\bu$ and the (half-turn) \emph{rotated} potential $Q^r$ (around the point $O$) by
\[\lb{ReQ'}
Q^r(v_{s,n_1,n_2})=
\left\{
\begin{array}{ll}
Q(v_{s,2-n_1,2-n_2}), & s=1, \\[2pt]
Q(v_{s,1-n_1,2-n_2}), & s=2, \\[2pt]
Q(v_{s,2-n_1,1-n_2}), & s=3,
\end{array}\right.
\]
for all $v\in\cV=\{v_{s,n_1,n_2}, \; s=1,2,3, \;n_1,n_2\in\Z\}$.

\begin{theorem}\lb{TKL} Let $H=A+Q$ be the Schr\"odinger operator with a $\G$-periodic complex potential $Q$ on the Kagome lattice $\bK$. Then, generically, the set of all $\G$-periodic potentials Floquet isospectral to $Q$ has the form
\[\lb{FIsK}
\mathrm{Iso}_\C^{\mathrm{F}}(Q)=\big\{Q_\ga,Q_\ga^r \mid \ga\in\G_o\},
\]
where $Q_\ga$ and $Q^r$ are the translated and rotated potentials defined by \er{TrQ'} and \er{ReQ'}, and $\G_o$ is given by \er{Ga_o'}.
\end{theorem}

\begin{remark}
The proofs of Theorems \ref{THL} and \ref{TKL} are based on the asymptotics \er{asy} and the Floquet spectral invariants $\cI_n^\mm(Q)$, $(n,\mm)\in\N\ts\Z^2$, for the discrete Schr\"odinger operators on periodic graphs obtained in \cite{S26}. Although the asymptotics \er{asy} holds for any periodic graphs, it is not a priori clear how many and what Floquet spectral invariants are enough to obtain the similar results for an arbitrary periodic graph. For example, for the hexagonal lattice (which is a bipartite graph) we have used only three invariants $\cI_{2p+2}^\mm(Q)$, $\mm=(1,0),(0,1),(-1,1)$, whereas for the Kagome lattice (which is a non-bipartite graph) we already need seven invariants $\cI_{2p+s}^{(1,0)}(Q)$, $\cI_{2p+s}^{(0,1)}(Q)$, $s=1,2,3$, and $\cI_{2p+2}^{(-1,1)}(Q)$.
\end{remark}

\section{Spectral invariants}\lb{Sec3}
\setcounter{equation}{0}
\subsection{Edge indices.}\lb{Sec2.1}
We define the important notion of an {\it edge index} which was introduced in \cite{KS14}. This notion allows one to consider, instead of a periodic graph, the finite quotient graph with edges labeled by some integer vectors called indices.

For each $\bx\in\R^d$, we denote by $\bx_\A\in\R^d$ the coordinate vector of $\bx$ with respect to the basis  $\A=\{\ga_1,\ldots,\ga_d\}$ of the lattice~$\G$, i.e.,
\[\lb{cola}
\bx_\A=(x_1,\ldots,x_d), \qqq \textrm{where} \qq \bx=
\textstyle\sum\limits_{s=1}^dx_s\ga_s.
\]

For any vertex $v\in\cV$ of a $\G$-periodic graph $\cG$, the following  unique representation holds true:
\[\lb{Dv}
v=v_0+[v], \qqq \textrm{where}\qqq v_0\in\cV\cap\Omega,\qqq [v]\in\G,
\]
and $\Omega$ is the fundamental cell of the lattice $\G$ defined by \er{fuce}. In other words, each vertex $v$ can be obtained from a vertex $v_0$ in $\Omega$  by a shift by a vector $[v]$ of the lattice $\G$.

For any oriented edge $\be=(u,v)\in\cA$ of the periodic graph $\cG$, we define the \emph{edge index} $\t(\be)$ as the vector of the lattice $\Z^d$ given by
$$
\t(\be)=[v]_\A-[u]_\A\in\Z^d,
$$
where $[v]\in\G$ is defined by \er{Dv} and the vector $[v]_\A\in\Z^d$  is given by \er{cola}. Due to the periodicity of the graph $\cG$, the edge indices $\t(\be)$ satisfy
\[\lb{Gpe}
\t(\be+\ga)=\t(\be),\qqq \forall\, (\be,\ga)\in\cA \ts\G.
\]
This periodicity of the indices allows us to assign the \emph{index} $\t(\be_*)\in\Z^d$ to each oriented edge $\be_*\in\cA_*$ of the quotient graph $\cG_*$ by setting
\[\lb{dco}
\t(\be_*)=\t(\be),
\]
where $\be\in\cA$ is an oriented edge in the periodic graph $\cG$ from the equivalence class $\be_*\in\cA_*=\cA/\G$. In other words, edge indices of the quotient graph $\cG_*$ are induced by edge indices of the periodic graph~$\cG$. Due to (\ref{Gpe}), the edge index $\t(\be_*)$ is uniquely determined by \er{dco} and does not depend on the choice of $\be\in\cA$. From the definition of the edge indices, it follows that
$$
\t(\ol\be\,)=-\t(\be), \qqq \forall\,\be\in\cA_*,
$$
where $\ol\be$ is the inverse edge of $\be\in\cA_*$.

\subsection{Closed paths and their characteristics} \lb{ssci}
A \emph{path} $\gp$ in a graph $\cG=(\cV,\cA)$ is a sequence of oriented edges
$$
\gp=(\be_1,\be_2,\ldots,\be_n), \qqq \textrm{where} \qqq
\be_j=(v_j,v_{j+1})\in\cA,\qq j=1,\ldots,n,
$$
for some vertices $v_1,\ldots,v_n,v_{n+1}\in\cV$. The number $n$ of edges in a path $\gp$ is called the \emph{length} of $\gp$ and is denoted by $|\gp|$, i.e., $|\gp|=n$. The \emph{reverse} of a path $\gp$ is the path $\ol\gp=(\ol\be_n,\ldots,\ol\be_1)$, where $\ol\be$ is the inverse edge of $\be\in\cA$. A path $\gp$ is called \emph{closed}, if $v_1=v_{n+1}$.

\begin{remark}\lb{Rcdv} If a graph $\cG$ has no multiple edges, then a closed path $\gp$ can equivalently be described by the sequence $\gp=(v_1,\ldots,v_n)$ of vertices it passes by.
\end{remark}

Let $\cC$ be the set of all closed paths in the quotient graph $\cG_*=(\cV_*,\cA_*)$. For any closed path $\bc\in\cC$ we define the \emph{index} $\t(\bc)\in\Z^d$ by
$$
\textstyle\t(\bc)=\sum\limits_{j=1}^n\t(\be_j),  \qqq \bc=(\be_1,\be_2,\ldots,\be_n)\in \cC,
$$
where $\t(\be)\in\Z^d$ is the index of the edge $\be\in\cA_*$.

At each vertex $v$ of the quotient graph $\cG_*=(\cV_*,\cA_*)$, we add a loop $\be_v$ with zero index and consider the \emph{modified} quotient graph
\[\lb{mfg}
\wt\cG_*=(\cV_*,\wt\cA_*),
\qqq \wt\cA_*=\cA_*\cup\{\be_v\}_{v\in\cV_*},\qqq \t(\be_v)=0.
\]

Let $\wt\cC$ be the set of all closed paths in $\wt\cG_*$, and $\ve>0$. For each closed path $\bc\in\wt\cC$, we define the \emph{weight}
\[\lb{Wcy}
\begin{array}{l}
\o(\bc,Q,\ve)=\o(\be_1)\ldots\o(\be_n), \qqq \textrm{where}\qqq \bc=(\be_1,\ldots,\be_n)\in\wt\cC,\\[6pt]
\o(\be)=\left\{
\begin{array}{cl}
\ve,  & \qq \textrm{if} \qq  \be\in\cA_* \\[3pt]
Q(v), & \qq \textrm{if} \qq \be=\be_v
\end{array}\right..
\end{array}
\]

\begin{remark}\lb{Rrep}
\emph{i}) The weight $\o(\bc,Q,\ve)$ of a closed path $\bc$ is a monomial of degree $|\bc|$ in the variables $\{\ve;\,Q(v),v\in\cV_*\}$.

\emph{ii}) The following identities hold true
$$
|\ol\bc|=|\bc|,\qqq \t(\ol\bc)=-\t(\bc), \qqq \o(\ol\bc,Q,\ve)=\o(\bc,Q,\ve), \qqq \forall\, \bc\in\wt\cC,
$$
where $\ol\bc$ is the reverse of $\bc$.
\end{remark}

Two closed paths are \emph{equivalent} if one is obtained by a cyclic permutation of edges in another. A \emph{cycle} is an equivalence class of a closed path $\bc=(\be_1,\ldots,\be_n)$ and will be denoted by $[\bc]=[\be_1,\ldots,\be_n]$. We define the length, index, and weight of $[\bc]$ as the length, index, and weight of any of its representatives.

\subsection{Spectral invariants for the rescaled Schr\"odinger operator}
Let $H=A+tQ$ be the Schr\"odinger operator with a periodic potential $Q$ on a periodic graph $\cG$, where $t>0$ is the coupling parameter. Instead of the operator $H$ it is convenient to introduce the rescaled Schr\"odinger operator
$$
H_\ve:=\frac1t\,H=\ve A+Q, \qqq \ve=\frac1t>0.
$$

\begin{lemma}\lb{Lres}
Periodic potentials $tQ$ and $tP$, $t>0$, on a fixed periodic graph $\cG$ are (Floquet) isospectral with respect to $A$ if and only if $Q$ and $P$ are (Floquet) isospectral with respect to~$\ve A$.
\end{lemma}
\no \textbf{Proof.} By the definition,
\begin{multline*}
\textrm{$tQ$ and $tP$ are Floquet isospectral with respect to $A$} \\ \Leftrightarrow\qq \s(A(k)+tQ)=\s(A(k)+tP)\qq \forall\, k\in\T^d\\
\Leftrightarrow\qq
\s(\ve A(k)+Q)=\s(\ve A(k)+P) \qq \forall\, k\in\T^d \\\Leftrightarrow
\qq\textrm{$Q$ and $P$ are Floquet isospectral with respect to $\ve A$.}\qqq \Box
\end{multline*}

\begin{proposition}\lb{PSpIn} Let $H_\ve=\ve A+Q$, $\ve>0$, be the rescaled Schr\"odinger operator with a periodic potential $Q$ on a periodic graph $\cG$ with the quotient graph $\cG_*=(\cV_*,\cE_*)$. Then\\
$\bu$ two potentials $Q_1$ and $Q_2$ are Floquet isospectral with respect to $\ve A$ if and only if
$$
\cI_n^\mm(Q_1,\ve)=\cI_n^\mm(Q_2,\ve),\qqq \forall\,(n,\mm)\in\N_\n\ts\Z^d,\qqq \n=\#\cV_*;
$$
\no $\bu$ two potentials $Q_1$ and $Q_2$ are simply isospectral with respect to $\ve A$ if and only if
$$
\cI_n(Q_1,\ve)=\cI_n(Q_2,\ve),\qqq \forall\,n\in\N_\n.
$$
Here
\[\lb{FspIn}
\cI_n^\mm(Q,\ve)=\frac1n\sum_{\bc\in\wt\cC_n^\mm\sm\cC_n^\mm}\o(\bc,Q,\ve),
\]
\[\lb{spIn}
\cI_n(Q,\ve)=\sum_{\mm\in\Z^d}\cI_n^\mm(Q,\ve)=\frac1n\sum_{\bc\in\wt\cC_n\sm \cC_n}\o(\bc,Q,\ve),
\]
where $\cC_n$ and $\wt\cC_n$ are the sets of all closed paths of length $n$ in the quotient graph $\cG_*$ and in the modified quotient graph $\wt\cG_*$ defined by (\ref{mfg}), respectively; \\
$\cC_n^\mm\ss\cC_n$ and $\wt\cC_n^\mm\ss\wt\cC_n$ are the sets of all closed paths of length $n$ and with index $\mm$ in $\cG_*$ and $\wt\cG_*$, respectively;\\
$\o(\bc,Q,\ve)$ is the weight of a closed path $\bc$ given by \er{Wcy}.
\end{proposition}

\no \textbf{Proof.} The proof is similar to the proof of Theorem 2.5 from  \cite{S26}, which is a particular case of this proposition for $\ve=1$. \qq $\Box$

\medskip

The functionals $\cI_n(Q,\ve)$ and $\cI_n^\mm(Q,\ve)$, $(n,\mm)\in\N_\n\ts\Z^d$, defined by \er{spIn} and \er{FspIn}, are called \emph{periodic} and \emph{Floquet} spectral invariants, respectively, for the rescaled Schr\"odinger operator~$H_\ve$.

\begin{proposition}\lb{Pfise}
i) The Floquet spectral invariants $\cI_n^\mm(Q,\ve)$, $(n,\mm)\in\N_\n\ts\Z^d$, defined by \er{FspIn} can be represented as finite power series in $\ve$:
\[\lb{Inme}
\cI_n^\mm(Q,\ve)=\sum_{k=\ell(\mm)}^{n-1}\ve^k\,\cI_{n,k}^\mm(Q),\qqq\textrm{where}\qqq
\cI_{n,k}^\mm(Q)=\frac1n\sum_{\bc\in\cC_k^\mm}W_{n-k}(\bc,Q),
\]
$\ell(\mm)$ is the length of the shortest closed path with index $\mm\in\Z^d$ in the quotient graph $\cG_*$; $\cC_k^\mm$ is the sets of closed paths $\bc=(v_1,\ldots,v_k)$ of length $k$ and with index $\mm$ in $\cG_*$, and $W_{n-k}(\bc,Q)$ is a homogeneous polynomial of degree $(n-k)$ in $Q(v_j)$, $j\in\N_k$.

In particular, for a closed path $\bc=(v_1,\ldots,v_k)$ with distinct vertices $v_1,\ldots,v_k$,
\[\lb{WnkQ}
W_{n-k}(\bc,Q)=\frac nk\,h_{n-k}(q_1,\ldots,q_k), \qqq q_j=Q(v_j),
\]
where $h_s(q_1,\ldots,q_k)$ is the complete homogeneous symmetric polynomial of degree $s$ in $q_1,\ldots,q_k$:
$$
h_s(q_1,\ldots,q_k)=\sum_{1\leq j_1\leq\ldots\leq j_s\leq k}
q_{j_1}\ldots q_{j_s}.
$$
Moreover, if each closed path in $\cC_k^\mm$ has distinct vertices, then
\[\lb{Incp}
\cI_{n,k}^\mm(Q)=\sum_{[v_1,\ldots,v_k]\in C_k^\mm}h_{n-k}(q_1,\ldots,q_k),
\]
where $C_k^\mm$ is the set of all cycles of length $k$ and with index $\mm$ in $\cG_*$.

ii) The periodic spectral invariants $\cI_n(Q,\ve)$, $n\in\N_\n$, defined by \er{spIn} can be represented as finite power series in $\ve$:
\[\lb{Ine}
\cI_n(Q,\ve)=\sum_{k=0}^{n-1}\ve^k\,\cI_{n,k}(Q),\qqq\textrm{where}\qqq
\cI_{n,k}(Q)=\frac1n\sum_{\bc\in\cC_k}W_{n-k}(\bc,Q),
\]
$\cC_k$ is the set of all closed paths of length $k$ in the quotient graph $\cG_*=(\cV_*,\cA_*)$. (A closed path of length 0 is any vertex of $\cG_*$.) In particular,
\[\lb{cIn0}
\cI_{n,0}(Q)=\frac1n\sum_{v\in\cV_*}Q^n(v);
\]
and, if $\cG_*$ has no loops, then
\[\lb{cIn12}
\cI_{n,1}(Q)=0,\qqq \cI_{n,2}(Q)=
\sum_{\{u,v\}\in\cE_*}m^2_{u,v}h_{n-2}(Q(u),Q(v)),
\]
where $\cE_*$ is the set of unoriented edges of $\cG_*$, and $m_{u,v}$ is the multiplicity of the edge $\{u,v\}\in\cE_*$.
\end{proposition}

\no \textbf{Proof.} \emph{i}) Due to \er{FspIn},
\[\lb{spIn1}
\cI_n^\mm(Q,\ve)=\frac1n\sum_{\bc\in\wt\cC_n^\mm\sm\cC_n^\mm}\o(\bc,Q,\ve),
\]
where $\o(\bc,Q,\ve)$ is the weight of $\bc$ given by \er{Wcy}.
Each closed path $\bc\in\wt\cC_n^\mm\sm\cC_n^\mm$ (of length $n$ and with index $\mm$) in the modified quotient graph $\wt\cG_*$ is a union of a closed path $\bc_o=(v_1,\ldots,v_k)\in\cC_k^\mm$ of length $\ell(\mm)\leq k\leq n-1$ and with index $\mm$ in the quotient graph $\cG_*$ and $n-k$ (not necessarily distinct) loops $\be_{v_{j_1}},\ldots,\be_{v_{j_{n-k}}}$, $j_1,\ldots,j_{n-k}\in\N_k$, added to the quotient graph $\cG_*$, see \er{mfg}. Recall that each such added loop $\be_v$ has the index $\t(\be_v)=0$ and the weight $\o(\be_v)=Q(v)$. Since each edge of $\cG_*$ has the weight $\ve$, then, due to \er{Wcy}, the weight of the closed path $\bc$ has the form
$$
\o(\bc,Q,\ve)=\ve^k\, q_{j_1}\ldots q_{j_{n-k}}, \qqq q_j=Q(v_j).
$$
Thus, using \er{spIn1}, we have
$$
\cI_n^\mm(Q,\ve)=\frac1n\sum_{k=\ell(\mm)}^{n-1}\sum_{\bc_o=(v_1,\ldots,v_k)\in\cC_k^\mm}
\ve^k\, W_{n-k}(\bc_o,Q),
$$
where
\[\lb{WnkQ'}
W_{n-k}(\bc_o,Q)=
\sum_{\bc=\bc_o\cup\{\be_{v_{j_1}},\ldots,\be_{v_{j_{n-k}}}\}\in\wt\cC_n^\mm}
q_{j_1}\ldots q_{j_{n-k}}
\]
is a homogeneous polynomial of degree $n-k$ in $q_j$, $j\in\N_k$. Thus, \er{Inme} has been proved.

If a closed path $\bc_o=(v_1,\ldots,v_k)\in\cC_k^\mm$ has distinct vertices $v_1,\ldots,v_k$, then
$$
W_{n-k}(\bc_o,Q)=\frac1k\sum_{\bc\in[\bc_o]}W_{n-k}(\bc,Q)=\frac nk
\sum_{1\leq j_1\leq\ldots\leq j_{n-k}\leq k}
q_{j_1}\ldots q_{j_{n-k}}=\frac nk\,h_{n-k}(q_1,\ldots,q_k).
$$
Here we have used that the cycle $[\bc_o]$ is an equivalence class consisting of $|\bc_o|$ distinct closed paths (these closed paths are a simple circular shift of one another). If each closed path in $\cC_k^\mm$ has distinct vertices, then, substituting \er{WnkQ} into the second identity in \er{Inme}, we get
$$
\cI_{n,k}^\mm(Q)=\frac1k\sum_{\bc=(v_1,\ldots,v_k)\in\cC_k^\mm}h_{n-k}(q_1,\ldots,q_k)=
\sum_{[\bc]=[v_1,\ldots,v_k]\in C_k^\mm}h_{n-k}(q_1,\ldots,q_k).
$$
Thus, \er{WnkQ}, \er{Incp} have also been proved.

\emph{ii}) The identity \er{Ine} is a direct consequence of the first identity in \er{spIn} and \er{Inme}. Now we prove \er{cIn0}, \er{cIn12}. A closed path of length 0 is a vertex of $\cG_*$. Then, using the second identity in \er{Ine} and \er{WnkQ'}, we have
$$
\cI_{n,0}(Q)=\frac1n\sum_{\bc\in\cC_0}W_n(\bc,Q)=\frac1n\sum_{v\in\cV_*}Q^n(v).
$$
If $\cG_*$ has no loops, then $\cC_1=\varnothing$, and, due to the second identity in \er{Ine},
$$
\cI_{n,1}(Q)=\frac1n\sum_{\bc\in\cC_1}W_{n-1}(\bc,Q)=0.
$$
Let $u,v\in\cV_*$, $u\neq v$, and let $\be_i=(u,v)\in\cA_*$, $i=1,\ldots,m_{u,v}$, where $m_{u,v}$ is the multiplicity of the edge $(u,v)$. Then each closed path $\bc=(u,v)$ of length 2 passing through the vertices $u,v$ has the form $(\be_i,\ol\be_j)$ for some (not necessarily distinct) $i,j=1,\ldots,m_{u,v}$. Thus, the number of such closed paths is $m_{u,v}^2$. Then, using the second identity in \er{Ine} and \er{WnkQ}, we obtain
\begin{multline*}
\cI_{n,2}(Q)=\frac1n\sum_{\bc\in\cC_2}W_{n-2}(\bc,Q)=
\frac12\sum_{\bc=(u,v)\in\cC_2}h_{n-2}(Q(u),Q(v))\\=
\frac12\sum_{(u,v)\in\cA_*}m_{u,v}^2h_{n-2}(Q(u),Q(v))=
\sum_{\{u,v\}\in\cE_*}m_{u,v}^2h_{n-2}(Q(u),Q(v)).
\end{multline*}
Thus, \er{cIn0}, \er{cIn12} have been proved. \qq $\Box$

\section{Proof of the main results}\lb{Sec4}
\setcounter{equation}{0}
\subsection{Isospectral potentials}
In this section we prove Proposition \ref{Pcar} about cardinality of the sets of isospectral potentials on periodic graphs and Proposition \ref{PAsy} about asymptotics for large isospectral potentials. In order to do this we need the following lemma, see \cite[Theorem 2.2]{F77}.

\begin{lemma}\lb{LF77}
Let $A$ be a given $\n\ts\n$ complex matrix. Then for any specified set $\s=\{\l_1,\ldots,\l_\n\}\ss\C$ there exists a diagonal complex matrix $D$ such that the spectrum of $A+D$ is the set $\s$. The number of such $D$ is finite and does not exceed $\n!$. Moreover, for almost all $\s$ the number of the diagonal matrices $D$ is exactly $\n!$.
\end{lemma}

\medskip

\no \textbf{Proof of Proposition \ref{Pcar}.} \emph{i}), \emph{ii})  For a given complex potential $Q$, the number $\#\Iso_\C(Q)$ is equal to the number of diagonal complex matrices $P$ such that $\s\big(A(0)+P\big)=\s\big(A(0)+Q\big)$, where $A(0)$ is the $\n\ts\n$ adjacency matrix of the quotient graph $\cG_*$. By Lemma~\ref{LF77}, the number of such $P$ does not exceed $\n!$. Moreover, for almost all $Q$ the number of $P$ is exactly $\n!$.

\emph{iii}) Let $Q=(q_1,\ldots,q_\n)$ be real and satisfy \er{sepQ}. Denote by $\l_1(Q)\leq\ldots\leq\l_\n(Q)$ the eigenvalues of $H(0)=A(0)+Q$. We rewrite the sequence $q_1,\ldots,q_\n$ in increasing order
$$
q_1^\bu<q_2^\bu<\ldots<q_\n^\bu.
$$
Here $q_1^\bu=q_{j_1}, q_2^\bu=q_{j_2},\ldots,q_\n^\bu=q_{j_\n}$
for some distinct integers $j_1,j_2,\ldots,j_\n\in\N_\n$. The eigenvalues of $H(0)$ satisfy
\[\lb{qq0}
q_j^\bu-\rho\le\l_j(Q)\le q_j^\bu+\rho, \qqq j\in\N_\n,
\]
where $\rho$ is the spectral radius of $A(0)$.

We fix a permutation $\pi$ of the set $\{1,\ldots,\n\}$ and consider the ball
\[\lb{Bro}
\cB_{\pi,\rho}=\big\{P=(p_1,\ldots,p_\n)\in\R^\n\mid \max_{j\in\N_\n}\big|p_{\pi(j)}-\l_j(Q)\big|\leq\rho\big\}.
\]
Let $P\in\cB_{\pi,\rho}$ and $i<j$. Then, using \er{Bro}, \er{qq0} and \er{sepQ}, we obtain
\begin{multline*}
p_{\pi(j)}-p_{\pi(i)}\geq\l_j(Q)-\rho-(\l_i(Q)+\rho)=\l_j(Q)-\l_i(Q)-2\rho\\\geq
q_j^\bu-\rho-(q_i^\bu+\rho)-2\rho=q_j^\bu-q_i^\bu-4\rho>0.
\end{multline*}
Thus, for $P=(p_1,\ldots,p_\n)\in\cB_{\pi,\rho}$,
\[\lb{incr}
p_{\pi(1)}<p_{\pi(2)}<\ldots<p_{\pi(\n)}.
\]
For every $P\in\cB_{\pi,\rho}$, the eigenvalues $\l_1(P)\leq\ldots\leq\l_\n(P)$ of $A(0)+P$ are real and satisfy
\[\lb{pp0}
p_{\pi(j)}-\rho\le\l_j(P)\le p_{\pi(j)}+\rho, \qqq j\in\N_\n.
\]

We define the continuous map $T:\R^\n\to\R^\n$ (depending on $\pi$) by
\[\lb{coT}
T(P)=(T_j(P))_{j=1}^\n, \qqq  T_{\pi(j)}(P)=p_{\pi(j)}-\l_j(P)+\l_j(Q),
\]
Using \er{coT} and \er{pp0}, for $P\in\cB_{\pi,\rho}$, we have
$$
|T_{\pi(j)}(P)-\l_j(Q)|=|p_{\pi(j)}-\l_j(P)|\leq\rho,\qqq \forall\,j\in\N_\n.
$$
Thus, $T(\cB_{\pi,\rho})\ss\cB_{\pi,\rho}$. Then, due to the Brouwer fixed point theorem, there exists a potential $P\in\cB_{\pi,\rho}$ such that $T(P)=P$, i.e.,
$$
p_{\pi(j)}-\l_j(P)+\l_j(Q)=p_{\pi(j)}, \qqq j\in\N_\n,
$$
which yields $\l_j(P)=\l_j(Q)$, $j\in\N_\n$, i.e., $P$ is isospectral to $Q$.

Due to the condition \er{incr}, the sets $\cB_{\pi,\rho}$ defined by \er{Bro} corresponding to distinct permutations $\pi$, do not intersect. Thus, there exist $\n!$ distinct potentials $P$ isospectral to $Q$. This and item \emph{i}) yield that there exist exactly $\n!$ isospectral potentials $P$ and each of $\n!$ regions defined by \er{repe} contains exactly one of them. \qq $\Box$

\medskip

\no\textbf{Proof of Proposition \ref{PAsy}.} Let $t>t_o$. Then
$$
\textstyle d(tQ)=td(Q)=t\,\frac{4\rho}{t_o}>4\rho.
$$
Thus, by Proposition \ref{Pcar}, the set of potentials isospectral to $tQ$ consists of $\n!$ real potentials. Let $tP$ be isospectral to $tQ$. Then, due to Lemma \ref{Lres}, $Q=(q_j)_{j=1}^\n$ and $P=(p_j)_{j=1}^\n$ are isospectral with respect to $\ve A$, $\ve=\frac1t$\,. According to Proposition \ref{PSpIn}, this is equivalent to the system of identities
\[\lb{IPIQ0}
\cI_n(P,\ve)=\cI_n(Q,\ve),\qqq \forall\,n\in\N_\n,
\]
where
\[\lb{InQe}
\cI_n(Q,\ve)=\sum_{k=0}^{n-1}\ve^k\,\cI_{n,k}(Q),
\]
and $\cI_{n,k}(Q)$ is a homogeneous polynomial of degree $(n-k)$ in $q_j$, $j\in\N_\nu$, see \er{Ine}. In particular, due to \er{cIn0},
\[\lb{InQ0}
\cI_n(Q,0)=\cI_{n,0}(Q)=\frac1n\sum_{j=1}^\n q_j^n.
\]

We consider the map $F=(F_1,\ldots,F_\n)$, $F_n=F_n(Q,\ve,P)$, where
\[\lb{FQeP}
F_n(Q,\ve,P)=\cI_n(P,\ve)-\cI_n(Q,\ve), \qqq n\in\N_\n.
\]
This map $F$ satisfies the conditions:

$\bu$ $F\in C^\iy(\R^{2\n+1})$,

$\bu$ $F(q_1,\ldots,q_\n,0,q_{\pi(1)},\ldots,q_{\pi(\n)})=0$ for each of $\n!$ permutations $\pi$ of the set $\{1,\ldots,\n\}$,

$\bu$ the Jacobian of the map $P\mapsto F(Q,0,P)$ is not equal to zero at each point $P=Q_\pi:=(q_{\pi(1)},\ldots,q_{\pi(\n)})$. Indeed, using \er{InQe} -- \er{FQeP},  we have
\begin{multline*}
\frac{\pa F_n}{\pa p_j}(Q,0,Q_\pi)=\frac{\pa\cI_n(P,\ve)}{\pa p_j}\bigg|_{\substack{P=Q_\pi\\\ve=0}}=\frac{\pa}{\pa p_j}\bigg(\sum_{k=0}^{n-1}\ve^k\,\cI_{n,k}(P)\bigg)\bigg|_{\substack{P=Q_\pi\\\ve=0}}=
\\=\frac{\pa\cI_{n,0}(P)}{\pa p_j}\bigg|_{P=Q_\pi}=\frac1n\,\frac{\pa}{\pa p_j}\bigg(\sum_{i=1}^\n p_i^n\bigg)\bigg|_{P=Q_\pi}=p_j^{n-1}\big|_{P=Q_\pi}=q_{\pi(j)}^{n-1},
\end{multline*}
which yields
\[\lb{Vade}
\bigg|\frac{\pa F}{\pa P}(Q,0,Q_\pi)\bigg|=\left|
\begin{array}{ccc}
1 & \ldots & 1 \\
q_{\pi(1)} & \ldots & q_{\pi(\n)} \\
q_{\pi(1)}^2 & \ldots & q_{\pi(\n)}^2 \\
\vdots & \ddots  & \vdots \\
q_{\pi(1)}^{\n-1} & \ldots & q_{\pi(\n)}^{\n-1}
\end{array}\right|=\prod_{1\leq i<j\leq\n}(q_{\pi(j)}-q_{\pi(i)})\neq0,
\]
since $q_i\neq q_j$, $i\neq j$. Then, due to the implicit function theorem, there exist an open set $U\ss\R^{\n+1}$ containing $(Q,0)$ and a unique system of functions $p_j:U\to\R$, $j\in\N_\n$, such that $P(Q,0)=Q_\pi$, and
$$
F(Q,\ve,P(Q,\ve))=0,\qqq \forall\,(Q,\ve)\in U.
$$
Moreover, $p_j$ are analytic functions. Thus, for $\ve$ small enough
\[\lb{ass0}
p_j=\sum_{r=0}^\iy p_j^{(r)}(Q)\,\ve^r,\qqq p_j^{(0)}(Q)=q_{\pi(j)},\qqq j\in\N_\n.
\]

We show that $p_j^{(1)}(Q)=0$ for all $j\in\N_\n$. Substituting \er{ass0} into the system \er{IPIQ0}, \er{InQe} and comparing the terms with $\ve^0$ and $\ve^1$, we obtain, using \er{cIn0} and the first identities in \er{cIn12},
$$
\sum_{j=1}^\n q_{\pi(j)}^n=\sum_{j=1}^\n q_j^n,\qqq \sum_{j=1}^\n q_{\pi(j)}^{n-1}p_j^{(1)}=0,
\qqq \forall\,n\in\N_\n.
$$
The first identities hold true. The second ones form the system of linear homogeneous equations with variables $p_j^{(1)}$, $j\in\N_\n$. Since its determinant, see \er{Vade}, is not zero, this system has the unique solution $p_j^{(1)}=0$, $j\in\N_\n$. Thus,
\[\lb{ass1}
p_j=q_{\pi(j)}+p_j^{(2)}\ve^2+O(\ve^3),\qqq j\in\N_\n.
\]
Substituting \er{ass1} into the system \er{IPIQ0}, \er{InQe} and comparing the terms with $\ve^2$, using \er{cIn0} and \er{cIn12}, we have
$$
\begin{array}{l}
\sum\limits_{j=1}^\n p_j^{(2)}=0,\qqq \sum\limits_{j=1}^\n q_{\pi(j)}p_j^{(2)}=0,\\[6pt]
\sum\limits_{j=1}^\n q_{\pi(j)}^{n-1}p_j^{(2)}+\cI_{n,2}(Q_\pi)=\cI_{n,2}(Q),\qq n=3,\ldots,\n, \qq Q_\pi=(q_{\pi(1)},\ldots,q_{\pi(\n)}),
\end{array}
$$
or, in the matrix form,
$$
V(Q_\pi)P^{(2)}=B(Q),
$$
where $V(Q_\pi)$ is the Vandermonde matrix, $P^{(2)}=\big(p_1^{(2)},\ldots,p_\n^{(2)}\big)^\top$,
$$
\begin{array}{ll}
B(Q)=\big(b_1(Q),\ldots,b_\n(Q)\big)^\top, & b_1(Q)=b_2(Q)=0,\\[8pt]
b_n(Q)=\cI_{n,2}(Q)-\cI_{n,2}(Q_\pi),\qqq & n=3,\ldots,\n.
\end{array}
$$
This yields
$$
P^{(2)}=V^{-1}(Q_\pi)B(Q),
$$
or, in the scalar form,
\[\lb{pi2}
p_i^{(2)}=\sum_{j=1}^\n\big(V^{-1}(Q_\pi)\big)_{ij}b_j(Q)=
\sum_{j=3}^\n\big(V^{-1}(Q_\pi)\big)_{ij}b_j(Q).
\]
Using the Vandermonde matrix inverse formula (see, e.g., \cite{MS58}), we have
\[\lb{VIij}
\big(V^{-1}(Q_\pi)\big)_{ij}=(-1)^{\n-j}\,
\frac{\gs_{\n-j}(q_{\pi(1)},\ldots,\hat q_{\pi(i)},\ldots,q_{\pi(\n)})}
{\prod\limits_{\substack{k=1\\ k\neq i}}^\n(q_{\pi(i)}-q_{\pi(k)})}\,,
\]
where $\gs_k(x_1,\ldots,x_n)$ is the elementary symmetric polynomial of degree $k$ in variables $x_1,\ldots,x_n$, and $\hat x_i$ means that the argument $x_i$ is excluded. Substituting \er{VIij} into \er{pi2} and using the second identity in \er{cIn12}, we obtain \er{coef2}. Thus, the asymptotics \er{asy} is proved. \qq $\Box$

\subsection{Floquet isospectrality on the hexagonal lattice}
In this section we consider the discrete Schr\"odinger operators with a $\G$-periodic potential $Q$ on the hexagonal lattice $\bG$. Recall that $\G$ is the lattice with the basis $p\ga_1,p\ga_2$ for some $p\in\N$, and $\ga_1,\ga_2$ are the periods of $\bG$, see Fig.~\ref{fig1}. We identify the vertices of the quotient graph $\bG_*=\bG/\G=(\cV_*,\cE_*)$ with the corresponding vertices of the hexagonal lattice $\bG$ in the fundamental cell of the period lattice~$\G$:
$$
\cV_*=\{v_{sij}, \; s=1,2; \; i,j\in\N_p\},\qqq q_{sij}:=Q(v_{sij}).
$$

\begin{proposition}\lb{PfsG}
Let $Q$ be a $\G$-periodic potential on the hexagonal lattice $\bG$. Then the length of the shortest cycles with indices $(1,0),(0,1),(-1,1)$ in the quotient graph $\bG_*=\bG/\G$ is $2p$, and the Floquet spectral invariants $\cI_{2p+2}^\mm(Q,\ve)$, $\mm=(1,0),(0,1),(-1,1)$, defined by \er{FspIn} have the form
\[\lb{ciFH}
\cI_{2p+2}^\mm(Q,\ve)=\ve^{2p}\,\cI_{2p+2,2p}^\mm(Q),
\]
where
\[\lb{F13H}
\begin{array}{l}
\cI_{2p+2,2p}^{(1,0)}(Q)=\sum\limits_{j=1}^p
h_2(q_{11j},q_{21j},q_{12j},q_{22j},\ldots,q_{1pj},q_{2pj})=:f_1(Q),\\[4pt] \cI_{2p+2,2p}^{(0,1)}(Q)=\sum\limits_{i=1}^p
h_2(q_{1i1},q_{2i1},q_{1i2},q_{2i2},\ldots,q_{1ip},q_{2ip})=:f_2(Q),\\[4pt]
\cI_{2p+2,2p}^{(-1,1)}(Q)=\sum\limits_{i=1}^p
h_2(q_{1i1},q_{2,i-1,1},q_{1,i-1,2},q_{2,i-2,2},\ldots,q_{1,i-p+1,p},q_{2,i-p,p})
=:f_3(Q),
\end{array}
\]
and $h_2(q_1,\ldots,q_k)$ is the complete homogeneous symmetric polynomial of degree 2 in $q_1,\ldots,q_k$.
\end{proposition}

\begin{remark}
Here and below the second and third indices $i,j$ in the vertices $v_{sij}$ of the quotient graph $\bG_*$ and in the corresponding values $q_{sij}:=Q(v_{sij})$ are taken modulo $p$.
\end{remark}

\no \textbf{Proof.} The length of the shortest cycles with indices $(1,0),(0,1),(-1,1)$ in $\bG_*$ is $2p$. All cycles of length $2p$ and with index $(1,0)$ in $\bG_*$ are the cycles
$$
\big[v_{11j},v_{21j},v_{12j},v_{22j},\ldots,v_{1pj},v_{2pj}\big],\qqq j\in\N_p
$$
(given by the sequences of their vertices). Similarly, all cycles of length $2p$ and with index $(0,1)$ in $\bG_*$ are
$$
\big[v_{1i1},v_{2i1},v_{1i2},v_{2i2},\ldots,v_{1ip},v_{2ip}\big],\qqq i\in\N_p.
$$
Finally, all cycles of length $2p$ and with index $(-1,1)$ in $\bG_*$ have the form
$$
\big[v_{1i1},v_{2,i-1,1},v_{1,i-1,2},v_{2,i-2,2},\ldots,v_{1,i-p+1,p},v_{2,i-p,p}\big],
\qqq i\in\N_p.
$$
The quotient graph $\bG_*$ is bipartite, then there are no cycles of odd length in $\bG_*$. Then, due to \er{Inme}, \er{Incp}, we obtain \er{ciFH}, \er{F13H}. \qq $\Box$

\medskip

We introduce the following vertex sets of the quotient graph $\bG_*=(\cV_*,\cE_*)$:

$\bu$ the $i$-\emph{th vertical}
\[\lb{cVvh}
\cV_\mathrm{ver}^i=\{v_{sij}\in\cV_*\mid s=1,2;\; j\in\N_p\},\qqq i\in\N_p,
\]
i.e., the set of vertices of the first ($s=1$) and second ($s=2$) types with the second index $i$;

$\bu$ the $j$-\emph{th horizontal}
\[\lb{cVhh}
\cV_\mathrm{hor}^j=\{v_{sij}\in\cV_*\mid s=1,2;\; i\in\N_p\},\qqq j\in\N_p,
\]
i.e., the set of vertices of the first and second types with the third index $j$;

$\bu$ the $m$-\emph{th diagonal}
\[\lb{cVdh}
\cV_\mathrm{diag}^m=\{v_{1ij},v_{2,i-1,j}\in\cV_*\mid i,j\in\N_p;\; i+j=m+1\},\qqq m\in\N_p,
\]
i.e., the set consisting of all vertices of the first and second types with the sum of the second and third indices equal to $m+1\,(\!\!\!\mod p)$ and $m\,(\!\!\!\mod p)$, respectively.

\medskip

The following lemma shows that a permutation $\pi:\cV_*\to\cV_*$ preserving the polynomials $f_s(Q)$, $s\in\N_3$, defined in \er{F13H} maps horizontals, verticals and diagonals onto horizontals, verticals and diagonals, respectively.

\begin{lemma}\lb{Lhvd'}
Let a permutation $\pi:\cV_*\to\cV_*$ preserve the polynomials $f_s(Q)$, $s\in\N_3$, defined in \er{F13H}. Then $\pi$ maps horizontals (verticals, diagonals) onto horizontals (verticals, diagonals, respectively), i.e., for any $i\in\N_p$ there exists $j\in\N_p$ such that
$$
\pi(\cV_\mathrm{hor}^i)=\cV_\mathrm{hor}^j \qqq \big(\pi(\cV_\mathrm{ver}^i)=\cV_\mathrm{ver}^j, \qq \pi(\cV_\mathrm{diag}^i)=\cV_\mathrm{diag}^j, \qq \textrm{respectively}\big).
$$
Here the sets $\cV_\mathrm{ver}^i$, $\cV_\mathrm{hor}^i$, $\cV_\mathrm{diag}^i$ are defined by \er{cVvh} -- \er{cVdh}.
\end{lemma}

\no \textbf{Proof.} Since the permutation $\pi:\cV_*\to\cV_*$ preserves the polynomial $f_1(Q)$, then
\begin{multline}\lb{ff33}
\textstyle\sum\limits_{j=1}^p
h_2\big(Q(\pi(v_{11j})),Q(\pi(v_{21j})),Q(\pi(v_{12j})),Q(\pi(v_{22j})),\ldots,
Q(\pi(v_{1pj})),Q(\pi(v_{2pj}))\big)\\[-6pt]\textstyle=
\sum\limits_{j=1}^ph_2(q_{11j},q_{21j},q_{12j},q_{22j},\ldots,q_{1pj},q_{2pj}).
\end{multline}
Let $j\in\N_p$ and $\pi(v_{11j})=v_{s'i'j'}$ for some $(s',i',j')\in\N_2\ts\N_p^2$. We write out all monomials in the LHS of \er{ff33} containing the factor $q_{s'i'j'}:=Q(v_{s'i'j'})$:
$$
q_{s'i'j'}
\big(Q(\pi(v_{11j}))+Q(\pi(v_{21j}))+Q(\pi(v_{12j}))+Q(\pi(v_{22j}))+\ldots+
Q(\pi(v_{1pj}))+Q(\pi(v_{2pj}))\big),
$$
and all monomials in the RHS of \er{ff33} containing the same factor $q_{s'i'j'}$:
$$
q_{s'i'j'}
\big(q_{11j'}+q_{21j'}+q_{12j'}+q_{22j'}+\ldots+q_{1pj'}+q_{2pj'}\big).
$$
Due to \er{ff33}, this yields
$$
\pi(v_{sij})\in\{v_{s'i'j'}\mid (s',i')\in\N_2\ts\N_p\},\qqq \forall\,(s,i)\in\N_2\ts\N_p,
$$
i.e., $\pi(\cV_\mathrm{hor}^j)=\cV_\mathrm{hor}^{j'}$. The remaining two identities follow from the invariance of the polynomials $f_2(Q),f_3(Q)$, respectively, under the permutation $\pi$. \qq $\Box$

\begin{proposition}\lb{Pv11'}
Let a permutation $\pi:\cV_*\to\cV_*$ preserve the polynomials $f_s(Q)$, $s\in\N_3$, defined in \er{F13H}, and $\pi(v_{111})=v_{111}$. Then $\pi$ is the identity map.
\end{proposition}

\no \textbf{Proof.} Since $\pi(v_{111})=v_{111}$ and, due to Lemma \ref{Lhvd'}, $\pi$ maps horizontals (verticals, diagonals) onto horizontals (verticals, diagonals, respectively), we have
\[\lb{1hvd}
\pi(\cV_\mathrm{hor}^1)=\cV_\mathrm{hor}^1,\qqq \pi(\cV_\mathrm{ver}^1)=\cV_\mathrm{ver}^1,\qqq \pi(\cV_\mathrm{diag}^1)=\cV_\mathrm{diag}^1.
\]
The proof consists of some steps.

\texttt{Step 1.} We show that $\pi$ preserves the vertices of the 1-st horizontal, i.e., $\pi(v_{si1})=v_{si1}$, $s=1,2$, $i\in\N_p$. The proof is by induction on $i$. Using \er{1hvd}, we have
$$
\pi(v_{211})\in\pi(\cV_\mathrm{hor}^1)\cap\pi(\cV_\mathrm{ver}^1)=
\cV_\mathrm{hor}^1\cap\cV_\mathrm{ver}^1=\{v_{111},v_{211}\}.
$$
Thus, $\pi(v_{211})=v_{211}$, and the base case $i=1$ holds true. Assume that $\pi(v_{sk1})=v_{sk1}$, $s=1,2$, for some $k<p$ and show that $\pi(v_{s,k+1,1})=v_{s,k+1,1}$, $s=1,2$. The induction hypothesis yields $\pi(\cV_\mathrm{diag}^{k+1})=\cV_\mathrm{diag}^{k+1}$. Then
$$
\pi(v_{1,k+1,1})\in\pi(\cV_\mathrm{hor}^1)\cap\pi(\cV_\mathrm{diag}^{k+1})=
\cV_\mathrm{hor}^1\cap\cV_\mathrm{diag}^{k+1}=\{v_{1,k+1,1},v_{2k1}\},
$$
which yields $\pi(v_{1,k+1,1})=v_{1,k+1,1}$ and $\pi(\cV_\mathrm{ver}^{k+1})=\cV_\mathrm{ver}^{k+1}$. Similarly,
$$
\pi(v_{2,k+1,1})\in\pi(\cV_\mathrm{hor}^1)\cap\pi(\cV_\mathrm{ver}^{k+1})=
\cV_\mathrm{hor}^1\cap\cV_\mathrm{ver}^{k+1}=\{v_{1,k+1,1},v_{2,k+1,1}\},
$$
and, consequently, $\pi(v_{2,k+1,1})=v_{2,k+1,1}$.

\texttt{Step 2.} Similarly to the proof of Step 1 we can show that $\pi$ preserves the vertices of the 1-st vertical, i.e., $\pi(v_{s1j})=v_{s1j}$, $s=1,2$, $j\in\N_p$.

\texttt{Step 3.} Since $\pi$ preserves the vertices of the 1-st horizontal and the 1-st vertical, then it preserves all horizontals, verticals and diagonals, i.e.,
\[\lb{pdvh}
\pi(\cV_\bu^i)=\cV_\bu^i,\qqq \bu=\mathrm{hor},\mathrm{ver},\mathrm{diag},\qqq i\in\N_p.
\]

\texttt{Step 4.} Finally, we show that $\pi$ is the identity map. Let $(i,j)\in\N^2_p$. Then, using \er{pdvh}, we obtain
$$
\begin{array}{l}
\pi(v_{1ij})\in\pi(\cV_\mathrm{ver}^i)\cap\pi(\cV_\mathrm{hor}^j)
\cap\pi(\cV_\mathrm{diag}^{i+j-1})=
\cV_\mathrm{ver}^i\cap\cV_\mathrm{hor}^j\cap\cV_\mathrm{diag}^{i+j-1}
=\{v_{1ij}\},\\[6pt]
\pi(v_{2ij})\in\pi(\cV_\mathrm{ver}^i)\cap\pi(\cV_\mathrm{hor}^j)
\cap\pi(\cV_\mathrm{diag}^{i+j})=
\cV_\mathrm{ver}^i\cap\cV_\mathrm{hor}^j\cap\cV_\mathrm{diag}^{i+j}
=\{v_{2ij}\}.
\end{array}
$$
Thus, $\pi(v_{sij})=v_{sij}$, $(s,i,j)\in\N_2\ts\N_p^2$. \qq $\Box$

\medskip

We call a permutation $\pi:\cV_*\to\cV_*$ the (half-turn) \emph{rotation}, if
\[\lb{syhe}
\pi(v_{1ij})=v_{2,p+1-i,p+1-j},\qqq \pi(v_{2ij})=v_{1,p+1-i,p+1-j},\qqq i,j\in\N_p.
\]

\begin{proposition}\lb{Pvpp}
Let a permutation $\pi:\cV_*\to\cV_*$ preserve the polynomials $f_s(Q)$, $s\in\N_3$, defined in \er{F13H}, and $\pi(v_{111})=v_{2pp}$. Then $\pi$ is the rotation defined by \er{syhe}.
\end{proposition}

\no \textbf{Proof.} The proof is similar to the proof of Proposition \ref{Pv11'}. Since $\pi(v_{111})=v_{2pp}$ and, due to Lemma \ref{Lhvd'}, $\pi$ maps horizontals (verticals, diagonals) onto horizontals (verticals, diagonals, respectively), we have
\[\lb{1hvd'}
\pi(\cV_\mathrm{hor}^1)=\cV_\mathrm{hor}^p,\qqq \pi(\cV_\mathrm{ver}^1)=\cV_\mathrm{ver}^p,\qqq \pi(\cV_\mathrm{diag}^1)=\cV_\mathrm{diag}^p.
\]

\texttt{Step 1.} We show that
\[\lb{cop1}
\pi(v_{1i1})=v_{2,p+1-i,p}\qqq\textrm{and}\qqq\pi(v_{2i1})=v_{1,p+1-i,p},\qqq i\in\N_p.
\]
The proof is by induction on $i$. Using \er{1hvd'}, we have
$$
\pi(v_{211})\in\pi(\cV_\mathrm{hor}^1)\cap\pi(\cV_\mathrm{ver}^1)=
\cV_\mathrm{hor}^p\cap\cV_\mathrm{ver}^p=\{v_{1pp},v_{2pp}\}.
$$
Thus, $\pi(v_{211})=v_{1pp}$, and the base case $i=1$ holds true. Assume that $$
\pi(v_{1k1})=v_{2,p+1-k,p}\qqq\textrm{and}\qqq\pi(v_{2k1})=v_{1,p+1-k,p}\qqq
\textrm{for some} \qq k<p,
$$
and show that
$$
\pi(v_{1,k+1,1})=v_{2,p-k,p}\qqq\textrm{and}\qqq\pi(v_{2,k+1,1})=v_{1,p-k,p}.
$$
The induction hypothesis yields $\pi(\cV_\mathrm{diag}^{k+1})=\cV_\mathrm{diag}^{p-k}$. Then
$$
\pi(v_{1,k+1,1})\in\pi(\cV_\mathrm{hor}^1)\cap\pi(\cV_\mathrm{diag}^{k+1})=
\cV_\mathrm{hor}^p\cap\cV_\mathrm{diag}^{p-k}=\{v_{1,p+1-k,p},v_{2,p-k,p}\},
$$
which yields $\pi(v_{1,k+1,1})=v_{2,p-k,p}$ (since $\pi(v_{2k1})=v_{1,p+1-k,p}$) and, consequently, $\pi(\cV_\mathrm{ver}^{k+1})=\cV_\mathrm{ver}^{p-k}$. Similarly,
$$
\pi(v_{2,k+1,1})\in\pi(\cV_\mathrm{hor}^1)\cap\pi(\cV_\mathrm{ver}^{k+1})=
\cV_\mathrm{hor}^p\cap\cV_\mathrm{ver}^{p-k}=\{v_{1,p-k,p},v_{2,p-k,p}\},
$$
which yields that $\pi(v_{2,k+1,1})=v_{1,p-k,p}$.

\texttt{Step 2.} Similarly, we can show that
\[\lb{cop2}
\pi(v_{11j})=v_{2,p,p+1-j}\qqq\textrm{and}\qqq\pi(v_{21j})=v_{1,p,p+1-j},\qqq j\in\N_p.
\]

\texttt{Step 3.} Since $\pi$ satisfies the conditions \er{cop1}, \er{cop2}, then
\[\lb{pdvh'}
\pi(\cV_\bu^i)=\cV_\bu^{p+1-i},\qqq \bu=\mathrm{hor},\mathrm{ver},\mathrm{diag},\qqq i\in\N_p.
\]

\texttt{Step 4.} Finally, we show that $\pi$ is given by \er{syhe}. Let $(i,j)\in\N^2_p$. Then, using \er{pdvh'}, we obtain
$$
\begin{array}{l}
\pi(v_{1ij})\in\pi(\cV_\mathrm{ver}^i)\cap\pi(\cV_\mathrm{hor}^j)
\cap\pi(\cV_\mathrm{diag}^{i+j-1})=
\cV_\mathrm{ver}^{p+1-i}\cap\cV_\mathrm{hor}^{p+1-j}\cap\cV_\mathrm{diag}^{p+2-i-j}
=\{v_{2,p+1-i,p+1-j}\},\\[6pt]
\pi(v_{2ij})\in\pi(\cV_\mathrm{ver}^i)\cap\pi(\cV_\mathrm{hor}^j)
\cap\pi(\cV_\mathrm{diag}^{i+j})=
\cV_\mathrm{ver}^{p+1-i}\cap\cV_\mathrm{hor}^{p+1-j}\cap\cV_\mathrm{diag}^{p+1-i-j}
=\{v_{1,p+1-i,p+1-j}\}.
\end{array}
$$
Thus, $\pi(v_{1ij})=v_{2,p+1-i,p+1-j}$ and $\pi(v_{2ij})=v_{1,p+1-i,p+1-j}$ for all $i,j\in\N_p$. \qq $\Box$

\medskip

We call a permutation $\pi_\ga:\cV_*\to\cV_*$ the \emph{translation} by a vector $\ga=l_1\ga_1+l_2\ga_2$, $l_j\in\N_p$, $j=1,2$, if
$$
\pi_\ga(v)=v+\ga \; (\hspace{-3mm}\mod \G), \qqq v\in\cV_*.
$$

Denote by $\cS$ the group of all permutations $\pi:\cV_*\to \cV_*$. Let $\cS_o\ss\cS$ be its subgroup generated by the translations $\pi_{\ga_1},\pi_{\ga_2}$ and the rotation defined by \er{syhe}.

\begin{proposition}\lb{PcSo'}
The subgroup $\cF$ of permutations in $\cS$, which leave the polynomials $f_s(Q)$, $s\in\N_3$, defined in \er{F13H}, invariant, is given by $\cS_o$.
\end{proposition}

\no \textbf{Proof.} First, we show that $\cS_o\subseteq\cF$. Let $\pi\in \cS_o$. Then $\pi$ is an automorphism of the quotient graph $\bG_*=(\cV_*,\cE_*)$ (i.e., an edge $\{u,v\}\in\cE_*$ if and only if $\{\pi(u),\pi(v)\}\in\cE_*$) preserving indices of all closed paths (the translations $\pi_\ga$) or changing their sign simultaneously (the rotation), i.e.,
$$
\t(\bc)=\t(\pi(\bc))\qq \forall\,\bc\in\cC,\qqq \textrm{or}\qqq
\t(\bc)=-\t(\pi(\bc))\qq \forall\,\bc\in\cC,
$$
where, for any closed path $\bc=(v_1,\ldots,v_n)$, the closed path $\pi(\bc)=(\pi(v_1),\ldots,\pi(v_n))$.

If $\pi$ preserves indices of all closed paths in $\bG_*$, then, $\bc\in\cC_k^\mm$ if and only if $\pi(\bc)\in\cC_k^\mm$. Recall that $\cC_k^\mm$ is the set of closed path of length $k$ and with index $\mm$ in $\bG_*$. Let the potential $Q_\pi$ be defined by $Q_\pi(v)=Q(\pi(v))$, $v\in\cV_*$. Then for each $(n,\mm)\in\N\ts\Z^2$ and $k<n$, due to the second identity in \er{Inme}, we have
$$
\cI_{n,k}^\mm(Q_\pi)=\frac1n\sum_{\bc\in\cC_k^\mm}W_{n-k}(\bc,Q_\pi)=
\frac1n\sum_{\pi(\bc)\in\cC_k^\mm}W_{n-k}(\pi(\bc),Q_\pi)=
\frac1n\sum_{\bc\in\cC_k^\mm}W_{n-k}(\bc,Q)=\cI_{n,k}^\mm(Q).
$$
Thus, the permutation $\pi$ preserves $f_s(Q)$, $s\in\N_3$, see \er{F13H}, and, consequently, $\pi\in\cF$.

Similarly, if $\pi\in\cS_o$ and $\pi$ changes the sign of indices of all closed paths in $\bG_*$ simultaneously, then $\bc\in\cC_k^\mm$ if and only if $\pi(\bc)\in\cC_k^{-\mm}$ and, consequently,
$$
\cI_{n,k}^\mm(Q_\pi)=\frac1n\sum_{\bc\in\cC_k^\mm}W_{n-k}(\bc,Q_\pi)=
\frac1n\sum_{\pi(\bc)\in\cC_k^{-\mm}}W_{n-k}(\pi(\bc),Q_\pi)=
\frac1n\sum_{\bc\in\cC_k^{-\mm}}W_{n-k}(\bc,Q)=\cI_{n,k}^{-\mm}(Q).
$$
For each closed path $\bc\in\cC_k^{-\mm}$, the reverse $\ol\bc\in\cC_k^{\mm}$ and $W_{n-k}(\ol\bc,Q)=W_{n-k}(\bc,Q)$, see Remark~\ref{Rrep}.\emph{ii}. This yields that $\cI_{n,k}^{-\mm}(Q)=\cI_{n,k}^{\mm}(Q)$. Thus, $\cI_{n,k}^\mm(Q_\pi)=\cI_{n,k}^\mm(Q)$, and $\pi\in\cF$.

Now we show that $\cF\subseteq\cS_o$. Let $\pi\in\cF$ and $\pi(v_{111})=v_{sij}$ for some $(s,i,j)\in\N_2\ts\N_p^2$. We consider the permutation
$$
\wt\pi=\pi_\ga\circ\pi, \qq\textrm{where}\qq
\ga=\left\{\begin{array}{ll}
(1-i)\ga_1+(1-j)\ga_2 \; (\hspace{-4mm}\mod \G), & \textrm{if } s=1,\\[4pt]
(p-i)\ga_1+(p-j)\ga_2 \; (\hspace{-4mm}\mod \G), & \textrm{if } s=2.
\end{array}\right.
$$
Then
$$
\wt\pi(v_{111})=\pi_\ga\big(\pi(v_{111})\big)=\pi_\ga(v_{sij})=
v_{sij}+\ga=\left\{\begin{array}{ll}
v_{111}, & \textrm{if } s=1,\\[4pt]
v_{2pp}, & \textrm{if } s=2.
\end{array}\right.
$$
Since $\pi\in\cF$ and $\pi_\ga\in\cS_o\subseteq\cF$, then $\wt\pi\in\cF$, and, due to Propositions \ref{Pv11'}, \ref{Pvpp}, $\wt\pi$ is either the identity map or the rotation defined by \er{syhe}. Since $\pi=(\pi_\ga)^{-1}\circ\wt\pi$ and $(\pi_\ga)^{-1},\wt\pi\in\cS_o$, then $\pi\in\cS_o$.\qqq $\Box$

\begin{remark}\lb{RigF}
From the proof of Proposition \ref{PcSo'} and the identities \er{Inme} it also follows that, for the hexagonal lattice, the invariance group of the Floquet spectral invariants $\cI_n^\mm(Q,\ve)$, $(n,\mm)\in\N_\n\ts\Z^2$, is given by $\cS_o$.
\end{remark}

\no \textbf{Proof of Theorem \ref{THL}.} Let periodic potentials $tQ$ and $tP$, $t>0$, be Floquet isospectral (with respect to $A$). Due to Lemma \ref{Lres}, $Q$ and $P$ are Floquet isospectral with respect to $\ve A$, $\ve=\frac1t$\,. Then, by Proposition \ref{PSpIn},
$$
\cI_{2p+2}^\mm(Q,\ve)=\cI_{2p+2}^\mm(P,\ve),\qqq \mm=(1,0),(0,1),(-1,1),
$$
or, using \er{ciFH}, \er{F13H},
\[\lb{f123}
f_s(Q)=f_s(P),\qqq s=1,2,3.
\]

Since $Q=(q_v)_{v\in\cV_*}$ and $P=(p_v)_{v\in\cV_*}$ are simply isospectral, then, due to Proposition \ref{PAsy}, for $\ve$ small enough
$$
p_v=q_{\pi(v)}+O(\ve^2),\qqq v\in\cV_*,
$$
for some permutation $\pi:\cV_*\to\cV_*$. Substituting this into \er{f123} and using \er{F13H}, we obtain
$$
f_s(Q)=f_s(Q_\pi), \qqq s\in\N_3, \qqq Q_\pi=(q_{\pi(v)})_{v\in\cV_*}.
$$
Thus, the permutation $\pi$ preserves the polynomials $f_s(Q)$, $s\in\N_3$, and, by Proposition \ref{PcSo'}, $\pi\in\cS_o$. Since each permutation $\pi$ from $\cS_o$ preserves the Floquet spectral invariants $\cI_n^\mm(Q,\ve)$, $(n,\mm)\in\N_\n\ts\Z^2$ (see Remark \ref{RigF}), then $P=Q_\pi$ and the set of all potentials Floquet isospectral to $Q$ has the form \er{FIsG}.\qq $\Box$

\subsection{Floquet isospectrality on the Kagome lattice}
In this section we consider the discrete Schr\"odinger operator with a $\G$-periodic potential $Q$ on the Kagome lattice $\bK$, where $\G$ is the lattice with the basis $p\ga_1,p\ga_2$ for some $p\in\N$, and $\ga_1,\ga_2$ are the periods of $\bK$, see Fig.~\ref{fig2}. We identify the vertices of the quotient graph $\bK_*=\bK/\G=(\cV_*,\cE_*)$ with the corresponding vertices of the Kagome lattice $\bK$ in the fundamental cell of the period lattice~$\G$:
$$
\cV_*=\{v_{sij}, \; s=1,2,3; \; i,j\in\N_p\},\qqq q_{sij}:=Q(v_{sij}).
$$

\begin{proposition}\lb{Pfsi}
Let $Q$ be a $\G$-periodic potential on the Kagome lattice $\bK$. Then the length of the shortest cycles with indices $(1,0),(0,1),(-1,1)$ in the quotient graph $\bK_*=\bK/\G$ is $2p$, and the Floquet spectral invariants $\cI_{2p+s}^{(1,0)}(Q,\ve)$, $\cI_{2p+s}^{(0,1)}(Q,\ve)$, $s=1,2,3$, and $\cI_{2p+2}^{(-1,1)}(Q,\ve)$, defined by \er{FspIn} have the form
\[\lb{ciFK}
\begin{array}{l}
\cI_{2p+1}^\mm(Q,\ve)=\ve^{2p}\,\cI_{2p+1,2p}^\mm(Q),\qqq \mm=(1,0),(0,1), \\[6pt]
\cI_{2p+2}^\mm(Q,\ve)=\ve^{2p}\,\cI_{2p+2,2p}^\mm(Q)+O(\ve^{2p+1}), \\[6pt]
\cI_{2p+3}^{\mm}(Q,\ve)=\ve^{2p}\,\cI_{2p+3,2p}^\mm(Q)+
\ve^{2p+1}\,\cI_{2p+3,2p+1}^\mm(Q)+O(\ve^{2p+2}),\\[6pt]
\cI_{2p+2}^{(-1,1)}(Q,\ve)=\ve^{2p}\,\cI_{2p+2,2p}^{(-1,1)}(Q)+O(\ve^{2p+1}),
\end{array}
\]
where
\[\lb{FI13}
\begin{array}{l}
\cI_{2p+s,2p}^{(1,0)}(Q)=
\sum\limits_{j=1}^ph_s(q_{11j},q_{21j},q_{12j},q_{22j},\ldots,q_{1pj},q_{2pj}),\\[4pt] \cI_{2p+s,2p}^{(0,1)}(Q)=
\sum\limits_{i=1}^ph_s(q_{1i1},q_{3i1},q_{1i2},q_{3i2},\ldots,q_{1ip},q_{3ip}),\qqq s=1,2,3,\\[4pt]
\cI_{2p+3,2p+1}^{(1,0)}(Q)=\sum\limits_{i,j=1}^p\big(h_2(q_{11j},q_{21j},\ldots,q_{1pj},q_{2pj},q_{3ij})+
h_2(q_{11j},q_{21j},\ldots,q_{1pj},q_{2pj},q_{3i,j-1})\big),\\[4pt]
\cI_{2p+3,2p+1}^{(0,1)}(Q)=\sum\limits_{i,j=1}^p
\big(h_2(q_{1i1},q_{3i1},\ldots,q_{1ip},q_{3ip},q_{2ij})+
h_2(q_{1i1},q_{3i1},\ldots,q_{1ip},q_{3ip},q_{2,i-1,j})\big),\\[4pt]
\cI_{2p+2,2p}^{(-1,1)}(Q)=\sum\limits_{i=1}^p
h_2(q_{2i1},q_{3i1},q_{2,i-1,2},q_{3,i-1,2},\ldots,q_{2,i-p+1,p},q_{3,i-p+1,p}),
\end{array}
\]
and $h_s(q_1,\ldots,q_k)$ is the complete homogeneous symmetric polynomial of degree $s$ in $q_1,\ldots,q_k$.
\end{proposition}

\no \textbf{Proof.} The length of the shortest cycles with indices $(1,0),(0,1),(-1,1)$ in the quotient graph $\bK_*$ is $2p$. Then, due to \er{Inme}, we obtain \er{ciFK}.

All cycles of length $2p$ and with index $(1,0)$ in $\bK_*$ are the cycles
\[\lb{cy10}
\bc_j:=[v_{11j},v_{21j},v_{12j},v_{22j},\ldots,v_{1pj},v_{2pj}],\qqq j\in\N_p
\]
(given by the sequences of their vertices). Then, due to \er{Incp},
$$
\cI_{2p+s,2p}^{(1,0)}(Q)=\sum_{j=1}^p
h_s(q_{11j},q_{21j},q_{12j},q_{22j},\ldots,q_{1pj},q_{2pj}),\qqq s=1,2,3.
$$

Each cycle of length $2p+1$ and with index $(1,0)$ in $\bK_*$ is obtained from one of the cycles $\bc_j$ given by \er{cy10} by replacing either an edge $(v_{1ij},v_{2ij})$, $i\in\N_p$, by two edges $(v_{1ij},v_{3ij})$ and $(v_{3ij},v_{2ij})$ or an edge $(v_{2ij},v_{1,i+1,j})$, $i\in\N_p$, by two edges $(v_{2ij},v_{3,i+1,j-1})$ and $(v_{3,i+1,j-1},v_{1,i+1,j})$. Thus, all cycles of length $2p+1$ and with index $(1,0)$ in $\bK_*$ have the form
\begin{multline*}
\big[v_{11j},v_{31j},v_{21j},v_{12j},v_{22j},\ldots,v_{1pj},v_{2pj}\big],\;
\big[v_{11j},v_{21j},v_{3,2,j-1},v_{12j},v_{22j},\ldots,v_{1pj},v_{2pj}\big],\;
\ldots\,,\\
\big[v_{11j},v_{21j},v_{12j},v_{22j},\ldots,v_{1pj},v_{3pj},v_{2pj}\big],\;
\big[v_{11j},v_{21j},v_{12j},v_{22j},\ldots,v_{1pj},v_{2pj},v_{3,1,j-1}\big],\qq j\in\N_p.
\end{multline*}
Then, using \er{Incp}, we get
$$
\cI_{2p+3,2p+1}^{(1,0)}(Q)=\sum\limits_{i,j=1}^p
\big(h_2(q_{11j},q_{21j},\ldots,q_{1pj},q_{2pj},q_{3ij})+
h_2(q_{11j},q_{21j},\ldots,q_{1pj},q_{2pj},q_{3,i,j-1})\big).
$$

The polynomials $\cI_{2p+s,2p}^{(0,1)}(Q)$, $s=1,2,3$, and $\cI_{2p+3,2p+1}^{(0,1)}(Q)$, $\cI_{2p+2,2p}^{(-1,1)}(Q)$ are derived in a similar way. We only note that
all cycles of length $2p$ and with index $(-1,1)$ in $\bK_*$ are the cycles
$$
\big[v_{2i1},v_{3i1},v_{2,i-1,2},v_{3,i-1,2},\ldots,v_{2,i-p+1,p},v_{3,i-p+1,p}\big],\qqq i\in\N_p.\qqq \Box
$$

\medskip

We introduce the following vertex set of the quotient graph $\bK_*=(\cV_*,\cE_*)$:

$\bu$ the set of $s$-\emph{type} vertices
\[\lb{cVs}
\cV_*^s=\{v_{sij}\in\cV_*\mid (i,j)\in\N_p^2\},\qqq s=1,2,3,
\]
i.e., the set of vertices with the first index $s$;

$\bu$ the $i$-\emph{th vertical}
\[\lb{cVv}
\cV_\mathrm{ver}^i=\{v_{sij}\in\cV_*\mid s=1,3,\; j\in\N_p\},\qqq i\in\N_p,
\]
i.e., the set of vertices of the first and third types with the second index $i$;

$\bu$ the $j$-\emph{th horizontal}
\[\lb{cVh}
\cV_\mathrm{hor}^j=\{v_{sij}\in\cV_*\mid s=1,2,\; i\in\N_p\},\qqq j\in\N_p,
\]
i.e., the set of vertices of the first and second types with the third index $j$;

$\bu$ the $m$-\emph{th diagonal}
\[\lb{cVd}
\cV_\mathrm{diag}^m=\{v_{sij}\in\cV_*\mid s=2,3;\; i,j\in\N_p;\; i+j=m+1\},\qqq m\in\N_p,
\]
i.e., the set of vertices of the second and third types with the sum of the second and third indices equal to $m+1\,(\!\!\!\mod p)$.

\medskip

For shortness, we introduce the notation
\[\lb{shno}
\begin{array}{lll}
f_1(Q):=\cI_{2p+1,2p}^{(1,0)}(Q),\qq & f_2(Q):=\cI_{2p+1,2p}^{(0,1)}(Q),\qq &\\[4pt]
f_3(Q):=\cI_{2p+2,2p}^{(1,0)}(Q), & f_4(Q):=\cI_{2p+2,2p}^{(0,1)}(Q),& f_5(Q):=\cI_{2p+2,2p}^{(-1,1)}(Q),\\[4pt]
f_6(Q):=\cI_{2p+3,2p+1}^{(1,0)}(Q),\qq & f_7(Q):=\cI_{2p+3,2p+1}^{(0,1)}(Q).&
\end{array}
\]
The following lemma shows that a permutation $\pi:\cV_*\to\cV_*$ preserving the polynomials $f_s(Q)$, $s\in\N_5$, preserves the vertex type and maps horizontals, verticals and diagonals onto horizontals, verticals and diagonals, respectively.

\begin{lemma}\lb{Lhvd}
Let a permutation $\pi:\cV_*\to\cV_*$ preserve the polynomials $f_s(Q)$, $s\in\N_5$, defined by \er{shno}, \er{FI13}. Then

i) the permutation $\pi$ preserves the vertex type, i.e.,
$$
\pi(\cV_*^s)=\cV_*^s, \qqq s=1,2,3;
$$

ii) the permutation $\pi$ maps verticals (horizontals, diagonals) onto verticals (horizontals, diagonals, respectively), i.e., for any $i\in\N_p$ there exists $j\in\N_p$ such that
$$
\pi(\cV_\mathrm{ver}^i)=\cV_\mathrm{ver}^j \qqq \big(\pi(\cV_\mathrm{hor}^i)=\cV_\mathrm{hor}^j, \qq \pi(\cV_\mathrm{diag}^i)=\cV_\mathrm{diag}^j, \qq \textrm{respectively}\big).
$$
Here the sets $\cV_*^s$, $\cV_\mathrm{ver}^i$, $\cV_\mathrm{hor}^i$, $\cV_\mathrm{diag}^i$ are given by \er{cVs} -- \er{cVd}.
\end{lemma}

\no \textbf{Proof.} \emph{i}) Since the permutation $\pi:\cV_*\to\cV_*$ preserves the polynomials $f_1(Q)$ and $f_2(Q)$, then
$$
\begin{array}{l}
\sum\limits_{i,j=1}^p\big(Q(\pi(v_{1ij}))+Q(\pi(v_{2ij}))\big)=
\sum\limits_{i,j=1}^p(q_{1ij}+q_{2ij}),\\
\sum\limits_{i,j=1}^p\big(Q(\pi(v_{1ij}))+Q(\pi(v_{3ij}))\big)=
\sum\limits_{i,j=1}^p(q_{1ij}+q_{3ij}).
\end{array}
$$
This yields that, for each $(i,j)\in\N_p^2$, $\pi(v_{1ij})\in\cV_*^1\cup\cV_*^2$ and $\pi(v_{1ij})\in\cV_*^1\cup\cV_*^3$, respectively. Thus,
$$
\pi(v_{1ij})\in(\cV_*^1\cup\cV_*^2)\cap(\cV_*^1\cup\cV_*^3)=\cV_*^1,\qqq \forall\,(i,j)\in\N_p^2,
$$
and, consequently, $\pi(v_{2ij})\in\cV_*^2$ and $\pi(v_{3ij})\in\cV_*^3$ for each $(i,j)\in\N_p^2$, i.e., $\pi$ preserves the vertex type.

\emph{ii}) Similarly, these statements follow from the invariance of the polynomials $f_3(Q),f_4(Q),f_5(Q)$, respectively, under the permutation $\pi$. \qq $\Box$

\medskip

We call the permutation $\pi:\cV_*\to\cV_*$ the (half-turn) \emph{rotation}, if
\[\lb{pesy}
\pi(v_{1ij})=v_{1,p+2-i,p+2-j},\;\; \pi(v_{2ij})=v_{2,p+1-i,p+2-j},\;\; \pi(v_{3ij})=v_{3,p+2-i,p+1-j},\;\; i,j\in\N_p.
\]

\begin{proposition}\lb{Pv11}
Let a permutation $\pi:\cV_*\to\cV_*$ preserve the polynomials $f_s(Q)$, $s\in\N_7$, defined by \er{shno}, \er{FI13}, and $\pi(v_{111})=v_{111}$. Then $\pi$ is either the identity map or the rotation given by \er{pesy}.
\end{proposition}

\no \textbf{Proof.} This is a direct consequence of Lemmas \ref{Lis1} -- \ref{Lis3}.\qq $\Box$

\begin{lemma}\lb{Lis1}
Let a permutation $\pi:\cV_*\to\cV_*$ preserve the polynomials $f_s(Q)$, $s\in\N_7$, defined by \er{shno}, \er{FI13}, and $\pi(v_{111})=v_{111}$. Then \[\lb{five}
\pi(v_{s11})=v_{s11},\; s=2,3, \qqq\textrm{or}\qqq
\pi(v_{211})=v_{2p1},\; \pi(v_{311})=v_{31p}.
\]
\end{lemma}

\no \textbf{Proof.} By Lemma \ref{Lhvd}, $\pi$ preserves the vertex type and maps horizontals (verticals, diagonals) onto horizontals (verticals, diagonals, respectively). Since $\pi(v_{111})=v_{111}$, then $\pi(\cV_\mathrm{hor}^1)=\cV_\mathrm{hor}^1$. This yields that
\[\lb{piv2}
\pi(v_{211})\in\{v_{2i1}\mid i\in\N_p\}.
\]
The permutation $\pi$ preserves the polynomial $f_7(Q)$ and $\pi(v_{111})=v_{111}$. Then $\pi$ preserves the sum of all monomials in $f_7(Q)$ containing the factor $q_{111}$:
$$
q_{111}\sum\limits_{j=1}^p
(2q_{111}+2q_{311}+2q_{112}+2q_{312}+\ldots+2q_{11p}+2q_{31p}+q_{21j}+q_{2pj}).
$$
From this it follows that
\[\lb{piv2'}
\pi(v_{211})\in\{v_{21j},v_{2pj}\mid j\in\N_p\}.
\]
Combining \er{piv2} and \er{piv2'}, we have
\[\lb{piv22}
\pi(v_{211})=v_{211} \qqq\textrm{or}\qqq
\pi(v_{211})=v_{2p1}.
\]
Similarly, using that $\pi(\cV_\mathrm{ver}^1)=\cV_\mathrm{ver}^1$ and $\pi$ preserves the polynomial $f_6(Q)$, we obtain
\[\lb{piv3}
\pi(v_{311})=v_{311} \qqq\textrm{or}\qqq
\pi(v_{311})=v_{31p}.
\]
Finally, since $\pi$ maps diagonals onto diagonals, from \er{piv22} and \er{piv3}, we derive \er{five}.\qq $\Box$

\begin{lemma}\lb{Lis2}
Let a permutation $\pi:\cV_*\to\cV_*$ preserve the polynomials $f_s(Q)$, $s\in\N_7$, defined by \er{shno}, \er{FI13}, and $\pi(v_{s11})=v_{s11}$, $s=1,2,3$. Then the permutation $\pi$ is the identity map.
\end{lemma}

\no \textbf{Proof.} Recall that, by Lemma \ref{Lhvd}, $\pi$ preserves the vertex type and maps horizontals (verticals, diagonals) onto horizontals (verticals, diagonals, respectively). The proof consists of some steps.

\texttt{Step 1.} We show that $\pi$ preserves the vertices of the 1-st horizontal, i.e., $\pi(v_{si1})=v_{si1}$, $s=1,2$, $i\in\N_p$. The proof is by induction on $i$. The base case $i=1$ holds true due to the condition of the lemma. Assume that $\pi(v_{sk1})=v_{sk1}$, $s=1,2$, for some $k<p$ and show that $\pi(v_{s,k+1,1})=v_{s,k+1,1}$, $s=1,2$. Since $\pi(v_{111})=v_{111}$, then $\pi(\cV_\mathrm{hor}^1)=\cV_\mathrm{hor}^1$. This yields that
\[\lb{v111}
\pi(v_{1,k+1,1})\in\{v_{1i1}\mid i\in\N_p\}.
\]
The permutation $\pi$ preserves the polynomial $f_7(Q)$ and $\pi(v_{2k1})=v_{2k1}$. Then $\pi$ preserves the sum of all monomials in $f_7(Q)$ containing the factor $q_{2k1}$:
$$
q_{2k1}
(q_{1k1}+q_{3k1}+\ldots+q_{1kp}+q_{3kp}+2q_{2k1}+
q_{1,k+1,1}+q_{3,k+1,1}+\ldots+q_{1,k+1,p}+q_{3,k+1,p}).
$$
From this it follows that
\[\lb{p1k1}
\pi(v_{1,k+1,1})\in\{v_{1kj},v_{1,k+1,j}\mid j\in\N_p\}.
\]
Combining \er{v111} and \er{p1k1} and using that, by the induction hypothesis, $\pi(v_{1k1})=v_{1k1}$, we obtain $\pi(v_{1,k+1,1})=v_{1,k+1,1}$.

Similarly, since $\pi(v_{1,k+1,1})=v_{1,k+1,1}$, $\pi$ preserves the sum of all monomials in $f_7(Q)$ containing the factor $q_{1,k+1,1}$:
$$
q_{1,k+1,1}\sum\limits_{j=1}^p
(2q_{1,k+1,1}+2q_{3,k+1,1}+\ldots+2q_{1,k+1,p}+2q_{3,k+1,p}+q_{2,k+1,j}+q_{2kj}),
$$
which yields that $\pi(v_{2,k+1,1})=v_{2,k+1,1}$.

\texttt{Step 2.} We can show that $\pi$ preserves the vertices of the 1-st vertical, i.e., $\pi(v_{s1j})=v_{s1j}$, $s=1,3$, $j\in\N_p$. The proof is similar to the proof of Step 1 and based on the invariance of the polynomial $f_6(Q)$ under the permutation $\pi$.

\texttt{Step 3.} Since $\pi$ preserves the vertices of the 1-st horizontal and the 1-st vertical, then it preserves all horizontals, verticals and diagonals, i.e.,
$$
\pi(\cV_\bu^i)=\cV_\bu^i,\qqq \bu=\mathrm{hor},\mathrm{ver},\mathrm{diag},\qqq i\in\N_p.
$$

\texttt{Step 4.} Finally, we show that $\pi$ is the identity map. Let $(i,j)\in\N^2_p$. Then the vertex $v_{1ij}$ belongs to the $i$-th vertical and $j$-th horizontal. Thus, $\pi(v_{1ij})=v_{1ij}$. Similarly, $\pi(v_{sij})=v_{sij}$, $s=2,3$. \qq $\Box$

\begin{lemma}\lb{Lis3}
Let a permutation $\pi:\cV_*\to\cV_*$ preserve the polynomials $f_s(Q)$, $s\in\N_7$, defined by \er{shno}, \er{FI13}, and
$$
\pi(v_{111})=v_{111}, \qqq \pi(v_{211})=v_{2p1},\qqq  \pi(v_{311})=v_{31p}
$$
Then the permutation $\pi$ is the rotation given by \er{pesy}.
\end{lemma}

\no \textbf{Proof.} The proof is similar to the proof of Lemma \ref{Lis2}. By Lemma \ref{Lhvd}, $\pi$ preserves the vertex type and maps horizontals (verticals, diagonals) onto horizontals (verticals, diagonals, respectively).

\texttt{Step 1.} We show that $\pi$ satisfies the conditions
\[\lb{cor1}
\pi(v_{1i1})=v_{1,p+2-i,1},\qqq \pi(v_{2i1})=v_{2,p+1-i,1},\qqq i\in\N_p.
\]
The proof is by induction on $i$. The base case $i=1$ holds true due to the condition of the lemma. Assume that $\pi(v_{1k1})=v_{1,p+2-k,1}$, $\pi(v_{2k1})=v_{2,p+1-k,1}$ for some $k<p$ and show that $\pi(v_{1,k+1,1})=v_{1,p+1-k,1}$, $\pi(v_{2,k+1,1})=v_{2,p-k,1}$. Since $\pi(v_{111})=v_{111}$, then $\pi(\cV_\mathrm{hor}^1)=\cV_\mathrm{hor}^1$. This yields that
\[\lb{v111a}
\pi(v_{s,k+1,1})\in\{v_{si1}\mid i\in\N_p\},\qqq s=1,2.
\]
We write out all monomials in $f_7(Q)$ containing the factor $q_{2,p+1-k,1}$:
\begin{multline*}
q_{2,p+1-k,1}
(q_{1,p+1-k,1}+q_{3,p+1-k,1}+\ldots+q_{1,p+1-k,p}+q_{3,p+1-k,p}\\+2q_{2,p+1-k,1}+
q_{1,p+2-k,1}+q_{3,p+2-k,1}+\ldots+q_{1,p+2-k,p}+q_{3,p+2-k,p}),
\end{multline*}
and all monomials in $f_7(Q_\pi)$, where $Q_\pi(v)=Q(\pi(v))$, $v\in\cV_*$, containing the same factor $q_{2,p+1-k,1}$:
\begin{multline*}
q_{2,p+1-k,1}
\big(q_{1,p+2-k,1}+Q(\pi(v_{3k1}))+\ldots+Q(\pi(v_{1kp}))+Q(\pi(v_{3kp}))\\
+2q_{2,p+1-k,1}+Q(\pi(v_{1,k+1,1}))+Q(\pi(v_{3,k+1,1}))+\ldots+Q(\pi(v_{1,k+1,p}))+
Q(\pi(v_{3,k+1,p}))\big).
\end{multline*}
Here we have used that $\pi(v_{1k1})=v_{1,p+2-k,1}$ and $\pi(v_{2k1})=v_{2,p+1-k,1}$. Since the permutation $\pi$ preserves the polynomial $f_7(Q)$, i.e., $f_7(Q)=f_7(Q_\pi)$, we obtain
\[\lb{p1k1a}
\pi(v_{1,k+1,1})\in\{v_{1,p+1-k,j},v_{1,p+2-k,j}\mid j\in\N_p\}.
\]
Combining \er{v111a} and \er{p1k1a} and using that, by the induction hypothesis, $\pi(v_{1k1})=v_{1,p+2-k,1}$, we obtain $\pi(v_{1,k+1,1})=v_{1,p+1-k,1}$.

Similarly, we write out all monomials in $f_7(Q)$ containing the factor $q_{1,p+1-k,1}$:
$$
q_{1,p+1-k,1}\sum\limits_{j=1}^p
(2q_{1,p+1-k,1}+2q_{3,p+1-k,1}+\ldots+2q_{1,p+1-k,p}+2q_{3,p+1-k,p}+q_{2,p+1-k,j}
+q_{2,p-k,j}),
$$
and all monomials in $f_7(Q_\pi)$ containing the same factor $q_{1,p+1-k,1}$:
\begin{multline*}
q_{1,p+1-k,1}\sum\limits_{j=1}^p
\big(2q_{1,p+1-k,1}+2Q(\pi(v_{3,k+1,1}))+\ldots+2Q(\pi(v_{1,k+1,p}))+
2Q(\pi(v_{3,k+1,p}))\\+Q(\pi(v_{2,k+1,j}))+Q(\pi(v_{2,k,j}))\big),
\end{multline*}
Here we have used that $\pi(v_{1,k+1,1})=v_{1,p+1-k,1}$. Since $f_7(Q)=f_7(Q_\pi)$, we obtain
\[\lb{p1k1b}
\pi(v_{2,k+1,1})\in\{v_{2,p+1-k,j},v_{2,p-k,j}\mid j\in\N_p\}.
\]
Combining \er{v111a} and \er{p1k1b} and using that, by the induction hypothesis, $\pi(v_{2k1})=v_{2,p+1-k,1}$, we obtain $\pi(v_{2,k+1,1})=v_{2,p-k,1}$.

\texttt{Step 2.} In a similar way, using the invariance of the polynomial $f_6(Q)$ under the permutation $\pi$, we can show that $\pi$ satisfies the conditions
\[\lb{cor2}
\pi(v_{11j})=v_{1,1,p+2-j},\qqq \pi(v_{31j})=v_{3,1,p+1-j},\qqq j\in\N_p.
\]

\texttt{Step 3.} Due to the first identities in \er{cor1}, we have
\[\lb{pive}
\pi(\cV_\mathrm{ver}^i)=\cV_\mathrm{ver}^{p+2-i},\qqq i\in\N_p.
\]
Similarly, from \er{cor2} we deduce that
\[\lb{pihd}
\pi(\cV_\mathrm{hor}^j)=\cV_\mathrm{hor}^{p+2-j}\qqq \textrm{and}\qqq
\pi(\cV_\mathrm{diag}^j)=\cV_\mathrm{diag}^{p+1-j},\qqq j\in\N_p.
\]

\texttt{Step 4.} Due to \er{pive} and \er{pihd}, we obtain that $\pi$ is the rotation defined by \er{pesy}. \qq $\Box$

\medskip

We call the permutation $\pi_\ga:\cV_*\to\cV_*$ the \emph{translation} by a vector $\ga=l_1\ga_1+l_2\ga_2$, $l_j\in\N_p$, $j=1,2$, if
$$
\pi_\ga(v)=v+\ga \; (\hspace{-3mm}\mod \G), \qqq v\in\cV_*.
$$

Denote by $\cS$ the group of all permutations $\pi:\cV_*\to \cV_*$. Let $\cS_o\ss\cS$ be its subgroup generated by the translations $\pi_{\ga_1},\pi_{\ga_2}$ and the rotation defined by \er{pesy}.

\begin{proposition}\lb{PcSo}
The subgroup $\cF$ of permutations in $\cS$, which leave the polynomials $f_s(Q)$, $s\in\N_7$, defined by \er{shno}, \er{FI13}, invariant, is given by $\cS_o$.
\end{proposition}

\no \textbf{Proof.} The proof of the inclusion $\cS_o\subseteq\cF$ repeats the arguments from the proof of Proposition~\ref{PcSo'}.

We prove that $\cF\subseteq\cS_o$. Let $\pi\in\cF$ and $\pi(v_{111})=v_{1ij}$ for some $(i,j)\in\N_p^2$. We consider the permutation
$$
\wt\pi=\pi_\ga\circ\pi, \qq\textrm{where}\qq \ga=(1-i)\ga_1+(1-j)\ga_2 \; (\hspace{-4mm}\mod \G).
$$
Then
$$
\wt\pi(v_{111})=\pi_\ga\big(\pi(v_{111})\big)=\pi_\ga(v_{1ij})=
v_{1ij}+(1-i)\ga_1+(1-j)\ga_2=v_{111}.
$$
Since $\pi\in\cF$ and $\pi_\ga\in\cF$, then $\wt\pi\in\cF$, and, due to Proposition \ref{Pv11}, $\wt\pi$ is either the identity map or the rotation defined by \er{pesy}. Since $\pi=(\pi_\ga)^{-1}\circ\wt\pi$ and $(\pi_\ga)^{-1},\wt\pi\in\cS_o$, then $\pi\in\cS_o$.\qq $\Box$

\medskip

\no \textbf{Proof of Theorem \ref{TKL}.} The proof is similar to the proof of Theorem \ref{THL}. Let periodic potentials $tQ$ and $tP$, $t>0$, be Floquet isospectral (with respect to $A$). Due to Lemma \ref{Lres}, $Q$ and $P$ are Floquet isospectral with respect to $\ve A$, $\ve=\frac1t$\,. Then, by Proposition \ref{PSpIn},
$$
\begin{array}{l}
\cI_{2p+s}^\mm(Q,\ve)=\cI_{2p+s}^\mm(P,\ve),\qqq \mm=(1,0),(0,1),\qqq  s=1,2,3,\\[6pt]
\cI_{2p+2}^{(-1,1)}(Q,\ve)=\cI_{2p+2}^{(-1,1)}(P,\ve),
\end{array}
$$
or, using \er{ciFK},
\[\lb{ciFK1}
\begin{array}{l}
\cI_{2p+1,2p}^\mm(Q)=\cI_{2p+1,2p}^\mm(P),\qqq \mm=(1,0),(0,1), \\[6pt]
\cI_{2p+2,2p}^\mm(Q)=\cI_{2p+2,2p}^\mm(P)+O(\ve), \\[6pt]
\cI_{2p+3,2p}^\mm(Q)+\ve\,\cI_{2p+3,2p+1}^\mm(Q)=\cI_{2p+3,2p}^\mm(P)+
\ve\,\cI_{2p+3,2p+1}^\mm(P)+O(\ve^2),\\[6pt]
\cI_{2p+2,2p}^{(-1,1)}(Q)=\cI_{2p+2,2p}^{(-1,1)}(P)+O(\ve),
\end{array}
\]
where the polynomials $\cI_{2p+s,2p}^\mm(Q)$, $s=1,2,3$, and  $\cI_{2p+3,2p+1}^\mm(Q)$, $\cI_{2p+2,2p}^{(-1,1)}(Q)$ are given in \er{FI13}.

Since $Q=(q_v)_{v\in\cV_*}$ and $P=(p_v)_{v\in\cV_*}$ are simply isospectral, then, due to Proposition \ref{PAsy}, for $\ve$ small enough
$$
p_v=q_{\pi(v)}+O(\ve^2),\qqq v\in\cV_*
$$
for some permutation $\pi:\cV_*\to\cV_*$. Substituting this into \er{ciFK1} and using \er{FI13}, we obtain
$$
\begin{array}{l}
\cI_{2p+1,2p}^\mm(Q)=\cI_{2p+1,2p}^\mm(Q_\pi)+O(\ve^2),\qqq \mm=(1,0),(0,1), \\[6pt]
\cI_{2p+2,2p}^\mm(Q)=\cI_{2p+2,2p}^\mm(Q_\pi)+O(\ve), \\[6pt]
\cI_{2p+3,2p}^\mm(Q)+\ve\,\cI_{2p+3,2p+1}^\mm(Q)=\cI_{2p+3,2p}^\mm(Q_\pi)+
\ve\,\cI_{2p+3,2p+1}^\mm(Q_\pi)+O(\ve^2),\\[6pt]
\cI_{2p+2,2p}^{(-1,1)}(Q)=\cI_{2p+2,2p}^{(-1,1)}(Q_\pi)+O(\ve),\qqq Q_\pi=(q_{\pi(v)})_{v\in\cV_*},
\end{array}
$$
which, in the notation \er{shno}, yields
$$
f_s(Q)=f_s(Q_\pi), \qqq s\in\N_7.
$$
Thus, the permutation $\pi$ preserves the polynomials $f_s$, $s\in\N_7$, and, by Proposition \ref{PcSo}, $\pi\in\cS_o$. Since each permutation from $\cS_o$ preserves the Floquet spectral invariants $\cI_n^\mm(Q,\ve)$, $(n,\mm)\in\N_\n\ts\Z^2$, then $P=Q_\pi$ and the set of all potentials Floquet isospectral to $Q$ has the form \er{FIsK}.\qq $\Box$

\medskip

\textbf{Acknowledgments.}  This work is supported by the Russian Science Foundation (project No. 25-21-00157).

\medskip

\textbf{Data availability}  No datasets were generated or analysed during the current study.

\medskip

\textbf{Conflict of interest} None declared.



\medskip

\begin{thebibliography}{9999}
\setlength{\itemsep}{-\parskip}\footnotesize
\bibitem[AM76]{AM76} Ashcroft N.W.; Mermin N.D. Solid State Physics, Holt, Rinehart and Winston, New York -- London, 1976.

\bibitem[BK13]{BK13} Berkolaiko, G.; Kuchment, P. Introduction  to
Quantum Graphs, Mathematical Surveys and Monographs, V. 186 AMS, 2013.

\bibitem[F77]{F77} Friedland, S. Inverse eigenvalue problems, Linear Algebra Appl., 17 (1977), no. 1, 15--51.

\bibitem[GKT93]{GKT93} Gieseker, D.; Kn\"orrer, H.; Trubowitz, E. The geometry of algebraic Fermi curves, volume 14 of Perspectives in Mathematics. Academic Press Inc., Boston, MA, 1993.

\bibitem[HN09]{HN09} Higuchi, Y.; Nomura, Y.
Spectral structure of the Laplacian on a covering graph,
European J. Combin. 30 (2009), no. 2, 570--585.

\bibitem[Ka88a]{Ka88a} Kappeler, T. On isospectral periodic potentials on a discrete lattice. I, Duke Math. J., 57 (1988), no. 1, 135--150.

\bibitem[Ka88b]{Ka88b} Kappeler, T. On isospectral potentials on a discrete lattice. II, Adv.in Appl. Math., 9 (1988), no. 4, 428--438.

\bibitem[Ka89]{Ka89} Kappeler, T. Isospectral potentials on a discrete lattice. III, Trans. Amer. Math. Soc., 314 (1989), no. 2, 815--824.

\bibitem[KS14]{KS14} Korotyaev, E.; Saburova, N.
Schr\"odinger operators on periodic discrete graphs, J. Math. Anal. Appl.,
 420 (2014), no. 1, 576--611.


\bibitem[KS22]{KS22} Korotyaev, E.; Saburova, N. Trace formulas for Schr\"odinger operators on periodic graphs, J. Math. Anal. Appl., 508 (2022), 33pp.

\bibitem[K16]{K16} Kuchment, P. An overview of periodic elliptic operators, Bull. Amer. Math. Soc., 53 (2016), 343--414.

\bibitem[L22]{L22} Liu, W. Topics on Fermi varieties of discrete periodic Schr\"odinger operators, J. Math. Phys., 63 (2022), 023503.

\bibitem[L23]{L23} Liu, W. Floquet isospectrality for periodic graph operators, J. Differ. Equ., 374 (2023), 642--653.

\bibitem[MS58]{MS58} Macon, N.; Spitzbart, A. Inverses of Vandermonde Matrices, Amer. Math. Monthly, 65 (1958), no. 2, 95--100.

\bibitem[RS78]{RS78} Reed, M.; Simon, B. Methods of modern mathematical
physics, vol.IV. Analysis of operators, Academic Press, New York, 1978.

\bibitem[S26]{S26} Saburova, N. Spectral invariants for discrete Schr\"odinger operators on periodic graphs, J. Spectr. Theory, 16 (2026), no. 2,  505--539.

\bibitem[SM25]{SM25} Saburova, N.; Martynov, D. Floquet isospectrality of the zero potential for Schr\"odinger operators on the hexagonal lattice, Proceedings of the International Conference Days on Diffraction, 2025, 188--193.

\bibitem[SS25]{SS25} Shipman, S.\,P.; Sottile, F. Algebraic aspects of periodic graph operators, 2025. arXiv:2502.03659.
\end{thebibliography}
\end{document}